\documentclass[11pt]{article}

\usepackage{amsmath,amssymb,latexsym}
\usepackage{theorem}
\usepackage{a4wide}
\usepackage[final]{graphicx}
\newcommand{\field}[1]{\mathbb{#1}}

\newcommand{\R}{\field{R}}
\newcommand{\Z}{\field{Z}}
\newcommand{\N}{\field{N}}

\newcommand{\DD}{{\mathcal D}}

\renewcommand{\S}{{\mathcal S}}

\newcommand{\card}{\mathop{\rm card}}

\newcommand{\isdef}{\stackrel{\text{\tiny def}}{=}}

\DeclareRobustCommand{\qed}{%
\ifmmode 
\else \leavevmode\unskip\penalty9999 \hbox{}\nobreak\hfill \fi
\quad\hbox{\qedsymbol}}
\newcommand{\openbox}{\leavevmode
\hbox to.77778em{%
\hfil\vrule
\vbox to.675em{\hrule width.6em\vfil\hrule}%
\vrule\hfil}}
\newcommand{\qedsymbol}{\openbox}
\newcommand{\proofname}{Proof}
\newenvironment{proof}[1][\proofname]{\par
\normalfont \trivlist \item[\hskip\labelsep   \itshape #1. ]
\ignorespaces
}{%
\qed\endtrivlist } 
\newtheorem{Theorem}{Theorem}
\newtheorem{Prop}[Theorem]{Proposition}

\newtheorem{corollary}[Theorem]{Corollary}
\newtheorem{Definition}[Theorem]{Definition}
{\theorembodyfont{\rmfamily}
\newtheorem{remark}{Remark}}
\begin{document}
\title{Discrete entropies of orthogonal polynomials}

\author{A.I.\ Aptekarev \and J.S.\ Dehesa \and A.
Mart\'{\i}nez-Finkelshtein\thanks{Corresponding author.} \and R.\ Y\'{a}\~{n}ez }

\maketitle
\smallskip

\begin{center}
\end{center}
\medskip

\begin{abstract}
Let $p_n$, $n\in \N$, be the $n$th orthonormal polynomial on $\R$, whose zeros are $\lambda _j^{(n)}$, $j=1, \dots, n$. Then for each $j=1, \dots, n$,
$$
\vec \Psi _j^2 \isdef  \left(\Psi_{1j}^2, \dots, \Psi_{nj}^2\right) 
$$
with
$$
\Psi _{ij}^2= p_{i-1}^2 (\lambda _j^{(n)}) \, \left( \sum_{k=0}^{n-1} p_k^2(\lambda_j^{(n)} )
\right)^{-1}    \,, \quad i=1, \dots, n,
$$
defines a discrete probability distribution. The Shannon entropy of the sequence $\{p_n\}$ is consequently defined as
$$
\mathcal S_{n,j}   \isdef   -\sum_{i=1}^n  \Psi _{ij}^{2} \log \left( \Psi
_{ij}^{2} \right)\,.
$$
In the case of Chebyshev polynomials of the first and second kinds an explicit and closed formula for $\mathcal S_{n,j}$ is obtained, revealing interesting connections with the number theory. Besides, several results of numerical computations exemplifying the behavior of $\mathcal S_{n,j}$ for other families are also presented.
\end{abstract}

\bigskip

\noindent AMS MOS Classification:\quad 33C45, 41A58, 42C05, 94A17
\medskip



\section{Introduction}
 
Given a probability Borel measure $\mu$ supported on the real line $\R$ with infinite number of points of increase, we can build a sequence of orthonormal polynomials $p_n(\lambda )=\kappa_n\, \lambda ^n +\text{lower degree terms}$, $n=0, 1, 2, \dots$, uniquely determined if all $\kappa _n>0$, such that
$$
\int p_n(\lambda ) p_m(\lambda )\, d\mu(\lambda )=\delta_{nm}\,, \quad m, n = 0, 1, 2, \dots
$$
Beside their importance in approximation theory and multiple branches of applied and pure mathematics, orthogonal polynomials constitute a noteworthy object from the point of view of the information theory. This interest originated in the framework of the modern density functional theory \cite{Dreizler90, Hohenberg64, March92}, that states that the physical and chemical
properties of fermionic systems (atoms, molecules, nuclei, solids)
may be completely described by means of the single-particle
probability density. For instance, if the solution of the
time-independent Schr\"{o}dinger equation in a $D$-dimensional
position space for an  single particle system,
$$
H \Psi (\vec{r}) = E \Psi(\vec{r})\,, \qquad \vec{r}=(x_1,
\dots, x_D)\,,
$$
is the wave function $\Psi(\vec{r})$,
then the position density of the system is $\rho(\vec{r}) = | \Psi(\vec{r}) |^2  $.
Analogously, the wave function in momentum space
$\widehat{\Psi}(\vec{p} )$, which is
the Fourier transform of
$\Psi(\vec{r} )$, gives the
momentum density $ \gamma(\vec{p}) = |
\widehat{\Psi}(\vec{p}) |^2 $.

Information measures of these densities are closely related to fundamental
and experimentally measurable physical quantities,  which makes them
useful in the study of the structure and dynamics of atomic and
molecular systems. For instance, the
\emph{Boltzmann-Gibbs-Shannon} (position-space) \emph{entropy}
\begin{equation} \label{eserho}
\mathcal B(\rho) = -\int \rho(\vec{r}) \log \rho(\vec{r}) \, d \vec{r}
\end{equation}
measures the uncertainty in the localization of  the particle in
space. Lower entropy corresponds to a more concentrated wave
function, with smaller uncertainty, and hence, higher
accuracy in predicting the localization of the particle. The well known inequality \cite{Babenko61, Bialynicki-Birula:84, Bialynicki-Birula:75}
\begin{equation} \label{xpeur}
\mathcal B(\rho) + \mathcal B(\gamma) \geq D(1+\log \pi)
\end{equation}
is an expression of the position-momentum uncertainty principle, much stronger than the renowned Heisenberg relation, that plays a major role in quantum mechanics (see \cite{Ohya93}).

The study of the information measures of orthogonal polynomials is motivated by the fact that the densities of many quantum mechanical systems with shape-invariant potentials (e.g., the harmonic oscillator and the hydrogenic systems) typically contain terms of the form $p_n^2 \, \mu'$. Explicit formulas, numerical algorithms and asymptotic behavior have been studied both for the Boltzmann-Gibbs-Shannon (or differential) entropy
$$
\mathcal B_n = -\int p_n^2 (\lambda ) \log \left( p_n^2(\lambda ) \mu '(\lambda )\right) \, d \mu (\lambda )
$$
and for the the relative entropy (or the Kullback-Leibler information)
$$
\mathcal K_n = -\int p_n^2 (\lambda ) \log \left( p_n^2(\lambda )  \right) \, d \mu (\lambda )\,;
$$
see e.g.\cite{Aptekarev:95, Aptekarev:96, Beckermann04, Buyarov04, MR2000f:94013, MR2001f:81036, MR1858268, MR2300963} and the references therein. In particular, it has been shown that for Chebyshev orthonormal polynomials of the first kind the relative entropy $\mathcal K_n$ does not depend on $n$,
$$
\mathcal K_n=\log ( 2)-1\,,
$$
and that this value is asymptotically maximal among all orthogonality measures on $[-1,1]$ (see \cite{Beckermann04}), giving a formal explanation to the intuitive notion that these polynomials are the most ``uniformly'' distributed ones.

However, there are several discrete measures naturally associated with a sequence of orthogonal polynomials. The analysis of such measures requires the use of the ``genuine'' entropy studied by Shannon. In order to stress the discrete character of this entropy, hereafter we refer to it as \emph{Shannon entropy} and denote it by the letter $\S$. The evaluation of Shannon entropy for discrete distributions is a basic question of information theory (see e.g.\ \cite{Jacquet99, Knessl98}); unlike for $\mathcal B_n$ and $\mathcal K_n$, there are no known results for the Shannon entropy of the orthogonal polynomials related distributions, due in part to the technical difficulties of the explicit evaluation of sums.

It is well known that the orthonormal polynomials $p_n$ satisfy a three-term recurrence relation of the form
\begin{equation}\label{recurrence}
\lambda \, p_i(\lambda )=b_{i+1}\, p_{i+1}(\lambda ) + a_{i+1}\,
p_i(\lambda ) + b_{i}\,  p_{i-1}(\lambda )\,, \quad i=0, 1, \dots,
n-2\,, \quad p_{-1}=0, \; p_0=1\,.
\end{equation}
Using its coefficients we can define the $n\times n$ Jacobi matrix ($n\in \N$),
\begin{equation}\label{JacobiMatrix}
   L_n=\begin{pmatrix} a_1 & b_1 &   \\
    b_1 & a_2 & b_1 &   \\
      & \ddots &\ddots & \ddots & \\
    &   & b_{n-2}  &  a_{n-1}  & b_{n-1} \\
&   &   &  b_{n-1} & a_{n}
\end{pmatrix}\,,
\end{equation}
which determines a self-adjoint linear operator (discrete Schr\"{o}dinger operator) $L_n:\, \R^n \to \R^n$
by
$$
L\,  \vec e_i = b_i \vec e_{i+1} + a_i \vec e_i + b_{i-1} \vec
e_{i-1}\,, \quad i=1, \dots, n\,,
$$
where $ \vec e_1, \dots \vec e_n$ is the canonical basis in $\R^n$, and we agree that $\vec e_0=\vec e_{n+1}=\vec
0$.
Moreover, up to a constant factor, $p_n(\lambda )=\det
(L_n -\lambda I)$, which shows that the eigenvalues $\lambda _k^{(n)}$, $k=1, \dots, n$,
are the zeros of $p_n$, and
$$
\vec P _k = \left(p_0, p_1(\lambda _k^{(n)}), \dots , p_{n-1}(\lambda _k^{(n)})\right)^T,
$$
are eigenvectors corresponding to different eigenvalues.

Let $\langle \cdot, \cdot \rangle$ denote the standard (euclidean) inner product in $\R^n$, and
\begin{equation}\label{def_ell}
    \ell_n(\lambda )\isdef \left( \sum_{k=0}^{n-1} p_k^2(\lambda )
\right)^{-1}=\frac{1}{\langle \vec P _k, \vec P _k \rangle }
\end{equation}
the $n$-th Christoffel function.

If we normalize
$$
\vec \Psi _k =\sqrt{\ell_n(\lambda_k^{(n)} )}\, \vec P _k\,,
$$
then a consequence of the well known Christoffel-Darboux formula is that
\begin{equation}\label{Mutualorthogonality}
     \langle \vec \Psi _i, \vec \Psi _j \rangle=\delta_{ij}\,,
\quad i,j =1,\dots, n\,.
\end{equation}
In other words, the $n\times n$ matrix
$$
\mathbf \Psi=\left( \sqrt{\ell_n(\lambda_j^{(n)} )}\,p_{i-1}(\lambda _j^{(n)})
\right)_{i,j=1}^n\,,
$$
made of columns $\vec \Psi _j$, $j=1, \dots, n$, is orthogonal, so that the squares of the components of each (column) vector $\vec \Psi _j$, $j=1, \dots, n$, give a discrete probability distribution, and these distributions are mutually orthogonal in the sense of \eqref{Mutualorthogonality}.

Recall that given a probability measure $\mu  = (\mu _1, \mu _2, \dots, \mu _n)$ on a system of $n$ points, e.\ g.\ $\sum_{j=1}^n \mu _j=1$, the standard \emph{Shannon entropy} reads $\mathcal S(\mu ) = - \sum_{j=1}^n  \mu _j \log \mu _j$.
By Jensen's inequality,
\begin{equation}
\label{ineq}
0\leq \mathcal S(\mu ) \leq \log (n)\,,
\end{equation}
and the maximum of corresponds to a uniform probability distribution. In this sense, it is quite
natural to think of the Shannon entropy as a measure of uncertainty.
 
\begin{remark}
We can give the following geometric interpretation to the Shannon entropy. Given in $\R^n$ an orthonormal basis $\{ \vec e_i\}$, any vector
$\vec v\in \R^n$ has a unique representation
$$
\vec v=\sum_{i=1}^n \langle \vec v, \vec e_i \rangle \, \vec e_i \,.
$$
Assume that $\vec v \in S^{n-1}$, that is, $\| \vec v\|=1$, where $\| \cdot \|$ means the Euclidean
norm. A natural way of measuring a relative distance of $\vec v$
from the basis $\{ \vec e_i\}$ is by means of the Shannon entropy
\begin{equation}\label{def_Entropy}
\mathcal S_n \isdef -\sum_{i=1}^n p_i \log (p_i)\,, \quad p_i\isdef \langle
\vec e_i ,\vec v  \rangle^2\,, \quad i=1, \dots, n\,.
\end{equation}
Indeed, if $\vec v=\vec e_k$ for a certain $k$, then
$p_j=\delta_{jk}$, and $\mathcal S_n=0$. On the contrary, if $\vec v$ is
``equidistant'' from all vectors $\vec e_j$'s, then all $p_j=1/n$,
and $\S_n$ attains its maximum, $\mathcal S_n=\log (n)$.
\end{remark}

Motivated by the discussion above, we introduce the \emph{discrete entropy of orthonormal polynomials} $p_n$, defined as the Shannon entropy of the probability distribution given by each column of $\mathbf \Psi$:
\begin{equation}\label{def_entropy_1}
\begin{split}
\mathcal S_{n,j} & \isdef   -\sum_{i=1}^n  \Psi _{ij}^{2} \log \left( \Psi
_{ij}^{2} \right) \\ &= - \ell_n(\lambda_j^{(n)} ) \, \sum_{i=1}^n
\,p_{i-1}^{2}(\lambda _j^{(n)}) \log \left( \ell_n(\lambda_j^{(n)} )
\,p_{i-1}^{2}(\lambda _j^{(n)}) \right) \\ &= - \log\left(\ell_n(\lambda_j^{(n)}
)\right) \, - \ell_n(\lambda_j^{(n)} ) \sum_{i=1}^n \,p_{i-1}^{2}(\lambda
_j^{(n)}) \log \left( p_{i-1}^{2}(\lambda _j^{(n)}) \right)  \,, \quad  j=1,
\dots, n\,,
\end{split}
\end{equation}
which can be generalized as
\begin{equation}\label{def_entropy_2}
\begin{split}
\mathcal  S_{n }(\lambda ) & \isdef   - \log\left(\ell_n(\lambda
)\right) - \ell_n(\lambda  ) \sum_{i=1}^n \,p_{i-1}^{2}(\lambda  )
\log \left( p_{i-1}^{2}(\lambda  ) \right)  \,,
\end{split}
\end{equation}
so that $\mathcal S_{n,j}=\mathcal S_{n}(\lambda _j^{(n)})$.

Unlike the Boltzmann-Gibbs-Shannon entropy, $\S_n$ does not depend on the weight function $\mu'$, and is suitable both for discrete and continuous orthogonality (cf.\ some numerical experiments in Section \ref{sec5}).

\begin{remark}
Since $\mathbf \Psi$ is an orthogonal matrix, its rows
$$
 \left( \sqrt{\ell_n(\lambda_1^{(n)} )}\,p_{i-1}(\lambda _1^{(n)}), \dots,
\sqrt{\ell_n(\lambda_n^{(n)} )}\,p_{i-1}(\lambda _n^{(n)}) \right)
$$
are also orthogonal vectors of $\R^n$:
\begin{equation*}\label{orthogonality2}
\delta_{ij}=  \sum_{k=1}^n \ell_n(\lambda_k^{(n)} ) \, p_{i-1}(\lambda
_k^{(n)})\, p_{j-1}(\lambda _k^{(n)}) =\int p_{i-1}(\lambda  ) p_{j-1}(\lambda
)\, d\mu_n(\lambda )\,,
\end{equation*}
where $\mu_n$ is the normalized counting measure of zeros of $p_n$:
\begin{equation*}\label{def_mu}
\mu_n = \sum_{k=1}^n  \ell_n(\lambda_k^{(n)} ) \delta_{\lambda _k^{(n)}}\,.
\end{equation*}
Hence, we may define the dual discrete entropy, corresponding to rows of $\mathbf \Psi$:
\begin{equation*}\label{def_entropy_3}
\begin{split}
\mathcal S_{n }^i & \isdef -  \sum_{j=1}^n
\ell_n(\lambda_j^{(n)} )  \,p_{i-1}^{2}(\lambda _j^{(n)}) \log \left(
\ell_n(\lambda_j^{(n)} ) \,p_{i-1}^{2}(\lambda _j^{(n)}) \right)   \,, \quad
j=1, \dots, n\,.
\end{split}
\end{equation*}

\end{remark}

\begin{figure}[htb]
\centering
\vspace{-2cm}\mbox{\includegraphics[scale=0.55]{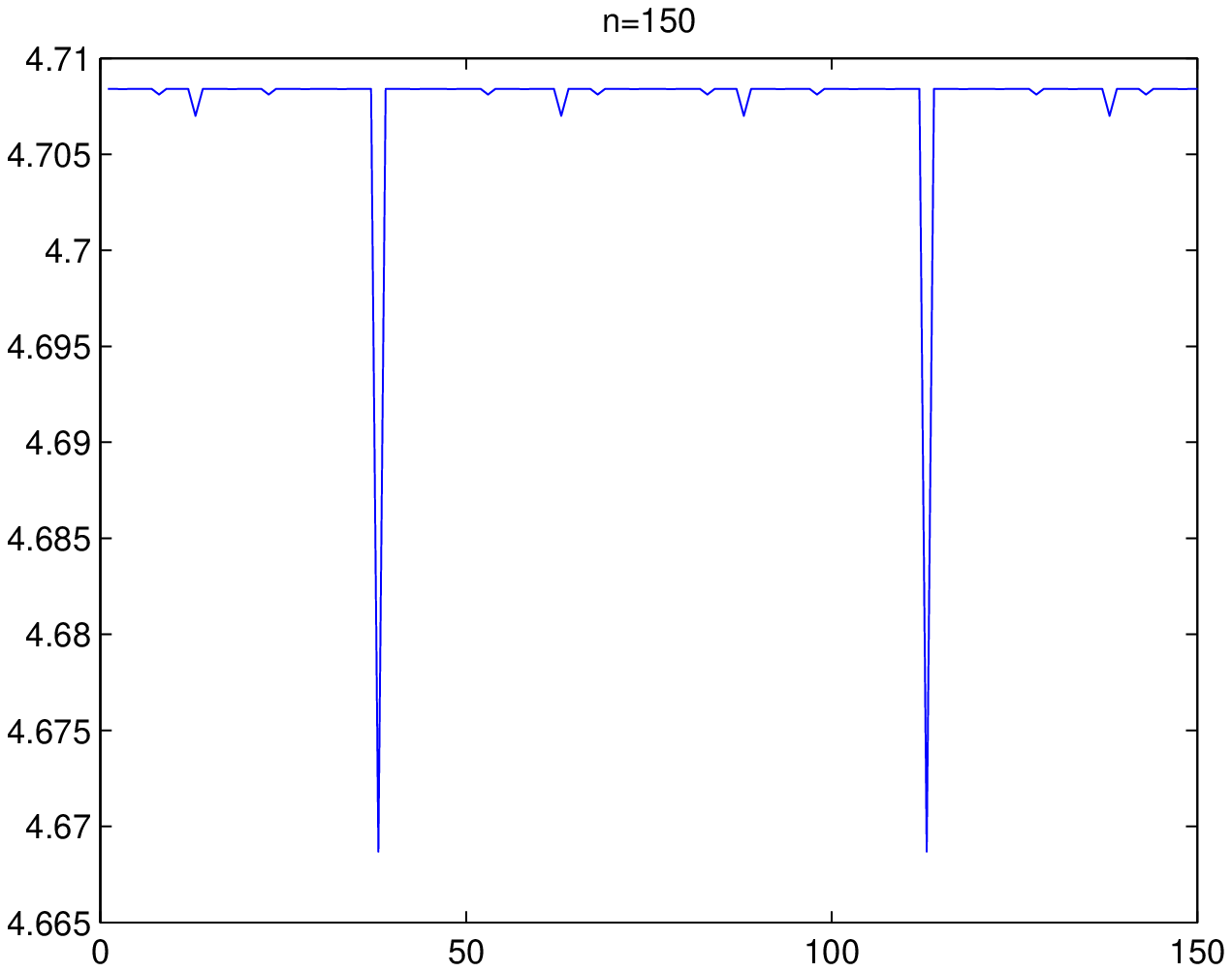}} \\
\vspace{-1cm}\mbox{\includegraphics[scale=0.55]{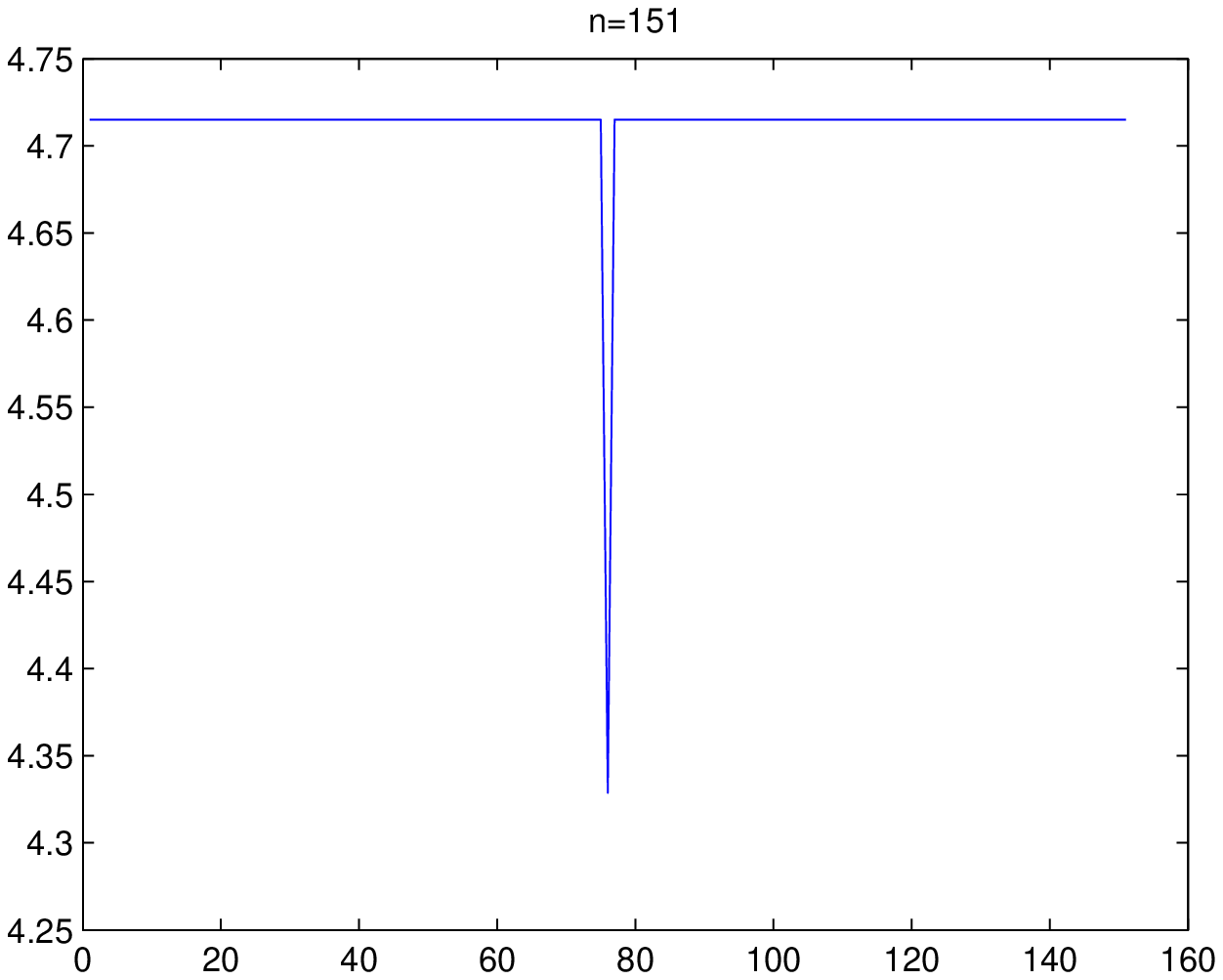}} \\
\vspace{-1cm}\mbox{\includegraphics[scale=0.55]{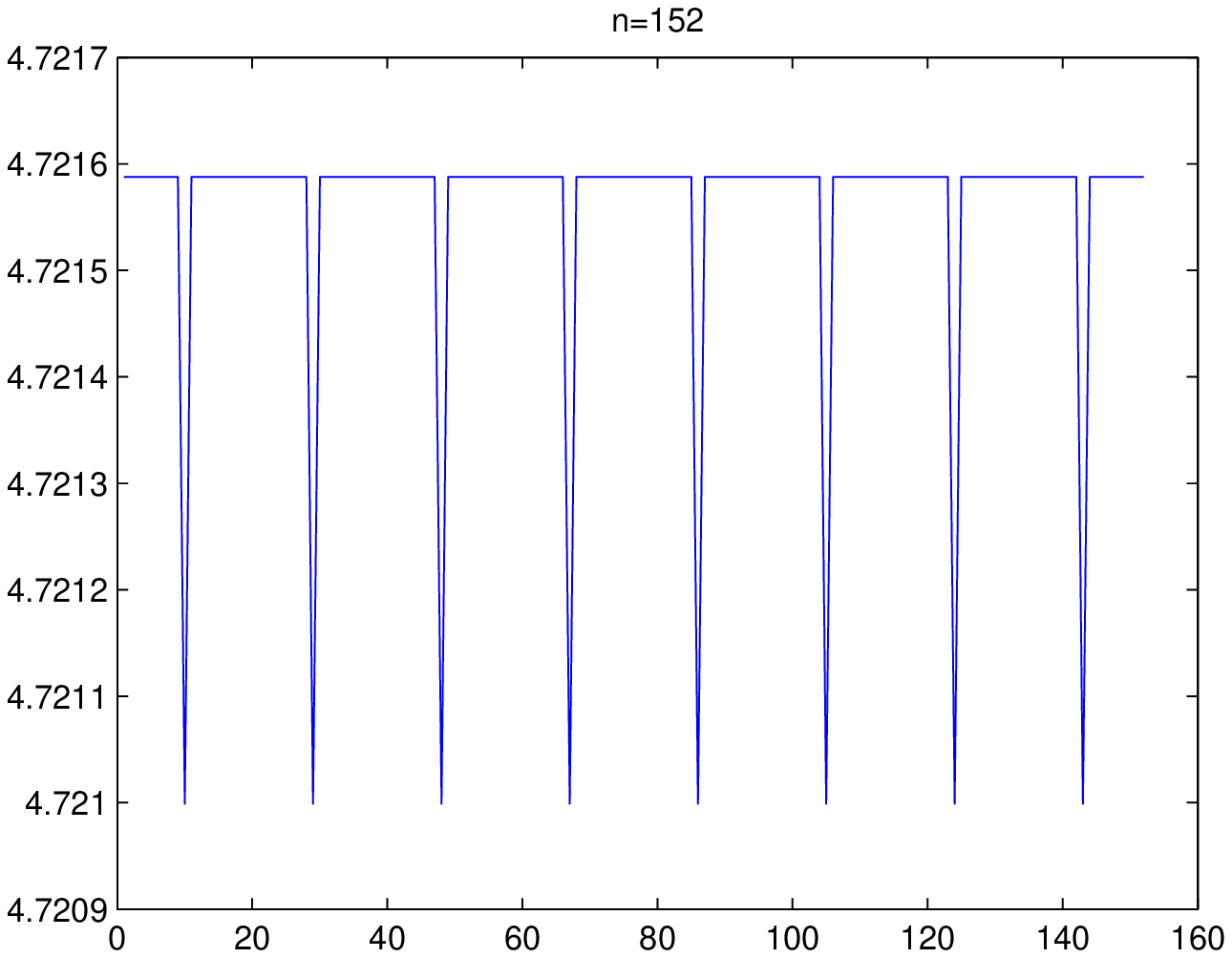}}
 \caption{Chebyshev polynomials of the first kind: entropy $S_{n,j}$ for $n=150$, $151$ and $152$.
}\label{fig:entropy150}
\end{figure}

A basic question of information theory is the evaluation of the Shannon entropy. In this paper we compute explicitly the discrete entropy $\mathcal S_{n,j}$ corresponding to Chebyshev orthonormal polynomials of the first and second kinds. A straightforward interpretation of \eqref{def_entropy_1} as Riemann sums allows to find the first two terms of the asymptotic expansion of $\mathcal S_{n,j}$ for fixed $j$ and large $n$; these terms do not depend on $j$. However, numerical experiments reveal the existence of certain picks, pointing downwards, whose position was not clear a priori (see Figure \ref{fig:entropy150}). The formulas presented below give a complete explanation of this phenomenon and exhibit nice connections with relevant objects from the number theory.

In order to state our results we need to introduce an auxiliary function
\begin{equation}\label{defR}
\mathcal R(x)\isdef x \left( \Psi\left(1-x\right) + 2 \gamma +
\Psi\left(1+x \right) \right)\,, \quad x\in [0,1)\,,
\end{equation}
where $\gamma $ is the Euler constant, and
$\Psi(x)=\Gamma'(x)/\Gamma(x)$ is the digamma function. Alternatively, $\mathcal R$ can be given by
its Taylor series expansion, absolutely convergent for $|x|<1$ (cf.\ formula (6.3.14) in
\cite{abramowitz/stegun:1972}),
\begin{equation}\label{series_for_R}
\mathcal R(x) = - 2 \sum_{k=1}^\infty \zeta(2k+1) x^{2k+1}\,, 
\end{equation}
where $\zeta(\cdot)$ is the Riemann zeta function.

Recall that \emph{Chebyshev polynomials of the first kind} are given by the explicit formula
\begin{equation}\label{explicitFirstKind}
   p_m(\lambda)= T_m(\lambda )=\begin{cases}  1, & \text{if } m=0\,, \\
\sqrt{2}\, \cos(m\theta), & \text{otherwise,}
\end{cases}\qquad \lambda=\cos\theta\,.
\end{equation}
They are orthonormal  with respect to the weight
$$
w(\lambda)=\frac{1}{\pi}\, \frac{1}{\sqrt{1-\lambda^2}} \quad
\text{on $[-1,1]$}\,.
$$
\begin{Theorem}
\label{corollary4} Let $n\in \N$, $j\in \{ 1, 2, \dots, n\}$. For
orthonormal Chebyshev polynomials of the first kind, the discrete
entropy has the following expression:
\begin{align*}
\mathcal S_{n,j} & = \log n + \log 2-1 +\frac{\log 2}{n} + \mathcal R
\left(\frac{d}{2 n}\right) \,, 
\qquad d=\text{GCD}(2j-1,n)\,.
\end{align*}
\end{Theorem}
Hereafter GCD stands for the greatest common divisor.

\begin{remark}
Observe that coefficients in the series expansion
\eqref{series_for_R} are all positive, so that  $\mathcal R(x)<0$
and is strictly decreasing for $x\in (0,1)$ (see Figure \ref{fig:psi}).
\begin{figure}[htb]
\centering
\hspace{-1.5cm}\mbox{\includegraphics[scale=0.9]{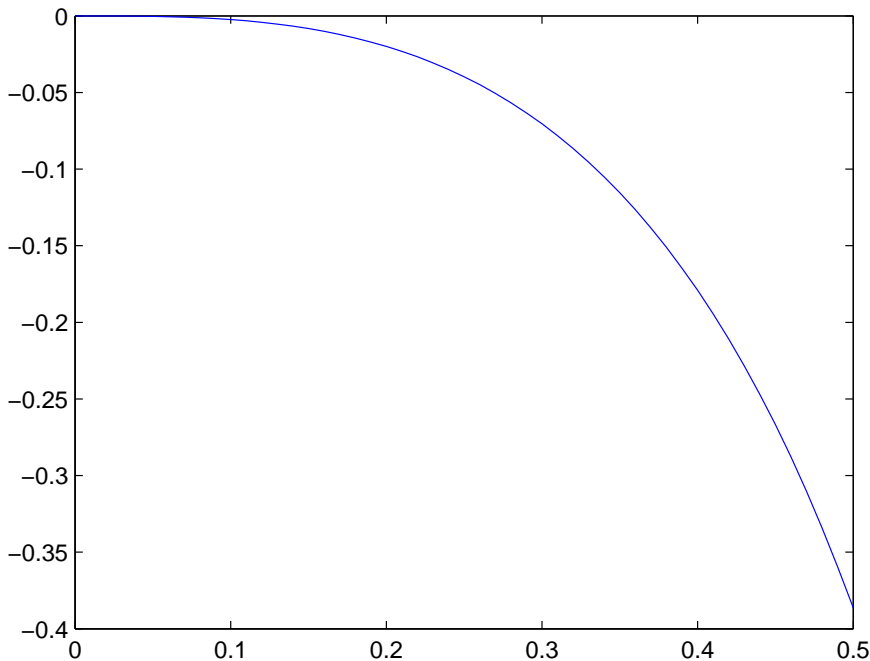}}
 \caption{Function $\mathcal R(x)$, for $x \in [0,1/2]$.
}\label{fig:psi}
\end{figure}
Since $1\leq d\leq n$, we see
that
\begin{align*}
\max_{j\in \{ 1, 2, \dots, n\}} \S_{n,j}=\log n + \log 2-1
+\frac{\log 2}{n} - 2 \sum_{k=1}^\infty \zeta(2k+1) \left(\frac{1}{2
n}\right)^{2k+1} \,,
\end{align*}
attained when $ \text{GCD}(2j-1,n)=1$. Furthermore, if $n$ is odd,
then $\S_{n,j}$ attains its minimum
\begin{align*}
\min_{j\in \{ 1, 2, \dots, n\}} \S_{n,j}=\log n -\log 2 +\frac{\log
2}{n} \,,
\end{align*}
at a single value $j=(n+1)/2$. It is the only local minimum of
$\S_{n,j}$ if $n\ge 3$ is prime.

The reader can compare these observations with the results of
numerical experiments shown in Figure \ref{fig:entropy150}.
\end{remark}

The \emph{Chebyshev polynomials of the second kind} are
\begin{equation}\label{explicitSecondKind}
  p_m(\lambda)=  U_m(\lambda )=\frac{\sin \left[ (m+1) \arccos (\lambda ) \right]}{\sqrt{1-\lambda ^2}} =
\frac{\sin \left[ (m+1) \theta \right]}{\sin (\theta )}\,, \quad \lambda=\cos\theta\,, \quad m \geq 0\,.
\end{equation}
They are orthonormal with respect to the weight
$$
w(\lambda )=\frac{2}{\pi}\,  \sqrt{1-\lambda ^2}  \quad
\text{on $[-1,1]$.}
$$
\begin{Theorem}
\label{corollary42Kind} Let $n\in \N$, $j\in \{ 1, 2, \dots, n\}$. For
orthonormal Chebyshev polynomials of the second kind, the discrete
entropy has the following expression:
\begin{align*}
\S_{n,j} & = \log (n+1) +  \log 2-1    + \mathcal R
\left(\frac{d}{  n+1}\right)\,, 
\qquad d=\text{GCD}(j,n+1)\,.
\end{align*}
\end{Theorem}

\begin{remark}
 
Since $1\leq d\leq n$, we see
that
\begin{align*}
\max_{j\in \{ 1, 2, \dots, n\}} \S_{n,j}=\log (n+1) + \log 2-1
 + \mathcal R
\left(\frac{1}{  n+1}\right) \,,
\end{align*}
attained when $ \text{GCD}(j,n+1)=1$. Furthermore, if $n$ is odd, $\S_{n,j}$ attains its minimum
\begin{align*}
\min_{j\in \{ 1, 2, \dots, n\}} \S_{n,j}=\log \left( \frac{n+1}{2}\right) \,,
\end{align*}
at a single value $j=(n+1)/2$. It is the only local minimum of
$\S_{n,j}$ if $n\ge 3$ is prime.
\end{remark}
\begin{remark}
The leading term $\log n$ in both cases shows that the values $p_0^2$, $p_1^2(\lambda _k^{(n)})$, \dots, $p_{n-1}^2(\lambda _k^{(n)})$, normalized by an appropriate factor, are approximately equidistributed. Comparing formulas from Theorem \ref{corollary4} and \ref{corollary42Kind} we see that unlike for the Boltzmann entropy, the discrete entropy of the Chebyshev polynomials of the first kind is generally smaller.
\end{remark}

The rest of the article is organized as follows. In the next section (that might have an independent interest) we discuss some piece-wise linear endomorphisms of $\R$ and their connection with permutations. This allows to reduce the analysis of the general discrete entropy $\S_{n,j}$ to some specific values of the index $j$. A modification of the Euler-Maclaurin summation formula is the key to the proof of Theorem \ref{corollary4} in Section \ref{subsection:1kind}. The close connection between polynomials of the first and second kinds allows us to avoid similar cumbersome computations in the proof of Theorem \ref{corollary42Kind} in Section \ref{subsection:2kind}. Finally, we discuss some numerical result obtained for the discrete entropy for other important families of orthogonal polynomials.

\section{Piece-wise linear endomorphism of $\R$ and permutations} \label{subsection:endomorphism}

\begin{figure}[htb]
\centering
\hspace{-1.5cm}\mbox{\includegraphics[scale=0.9]{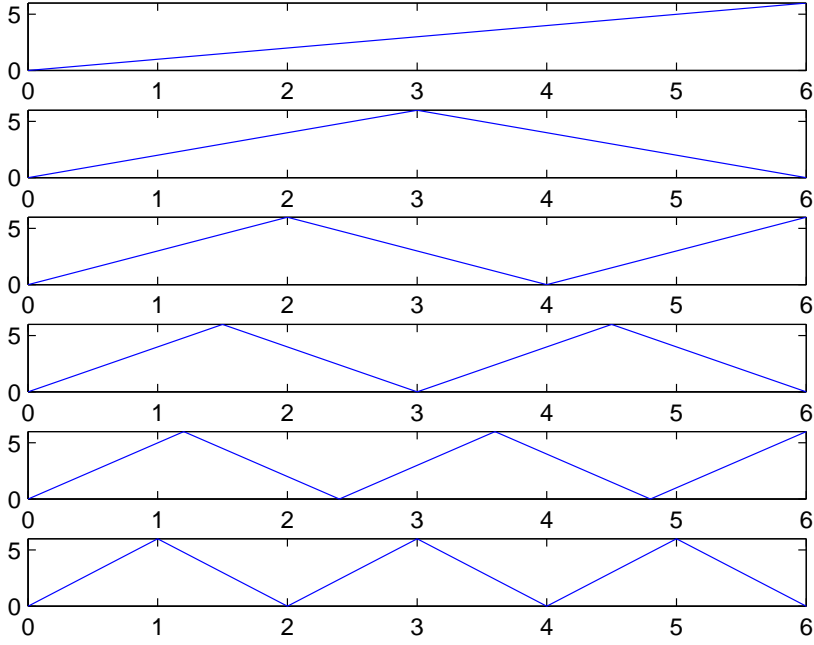}}
 \caption{Functions $\varphi_j^{(6)}(x)$, $x\in [0,6]$, for $j=1, \dots, 6$.
}\label{fig:phi1}
\end{figure}

The key role is played by the following auxiliary function:
\begin{Definition}
For each pair of values $n, j\in \N$,   
let $\varphi_j^{(n)}(x)
$ denote the linear spline on $\R$ with nodes at $\{ mn/j\}_{m\in \Z}$ interpolating the values
$$
\left \{ \left(\frac{ m n}{ j }, n \, \frac{1-(-1)^m }{2} \right)\right\}_{m\in \Z}\,.
$$
\end{Definition}
Functions $\varphi_j^{(n)}$ can be explicitly described by
\begin{equation}\label{alt1}
x\in\left[ \dfrac{2 k-1 }{j}\,
n,\dfrac{2 k+1}{j}\, n\right] \text{ for } k \in \Z \quad \Rightarrow \quad \varphi_j^{(n)}(x ) =  |jx-  2 k n|
\end{equation}
(see Figure \ref{fig:phi1}). 

In order to summarize necessary properties of functions $\varphi_j^{(n)}$ we need to introduce some notation.
We denote by $a\equiv b \mod (c)$ the standard arithmetic congruence of $a$ and $b$ modulo $c$, $ \text{GCD}(a,b) $ stands for the greatest common divisor of integer numbers $a$ and $b$, and  $\N_0 \isdef \N\cup \{ 0\}$. We define also the remainder function
$\mathcal D:\, \N_0 \times \N  \to \Z$, by
$$
\DD (p,q)= r \quad \text{if and only if} \quad -q/2 \leq r < q/2 \text{ and
} p \equiv r \mod(q)
$$
(note that this definition is shifted with respect to a standard concept of remainder; thus, $\DD (p,q)$ takes also negative values).

\medskip

\noindent \textbf{Main Lemma} \emph{Let $n, j\in \N$,  with $\text{GCD}(j,n)=d$.}
\begin{enumerate}
\item[(i)] \emph{We have
\begin{equation}\label{symmetryPhiMore}
    \varphi_j^{(n)}\left( x+\frac{n}{d}\right) = \begin{cases}
     n- \varphi_j^{(n)}\left( x\right)\,, & \text{if $j/d$ is odd,}\\
     \varphi_j^{(n)}\left( x\right)\,, & \text{if $j/d$ is even,}
    \end{cases}
    \quad x \in \R\,,
\end{equation}
and for $k\in \Z$,
\begin{equation}\label{hiddenSymbis}
\varphi_j^{(n)}(k )=d\cdot \varphi_{j/d}^{(n/d)}(k )\,.
\end{equation}}

\item[(ii)] \emph{If $d=n$, then $\varphi_j^{(n)}(\Z)\subset \{ 0, n\}$.}

\item[(iii)] \emph{If $\text{GCD}(j,2n)=d<n$, then
\begin{equation}\label{identity1bis}
 \varphi_j^{(n)} \left( \left\{ 1,\dots, n-1\right\}\right)\setminus \{ 0, n\} = \left\{d m:\, m=1 ,\dots,
\frac{n}{d}-1\right\}\,,
\end{equation}
and for any
$m\in\left\{1,\dots,\frac{n}{d}-1\right\} $,
\begin{equation}\label{identity2bis}
\card \{k \in\{1,\ldots,n-1\}:\varphi_j^{(n)}(k)=d m\}=d\,.
\end{equation}}
\item[(iv)] \emph{If $\text{GCD}(j,2n)=2 d$ and $d<n$, then
\begin{equation}\label{identity1even}
 \varphi_j^{(n)} \left( \left\{ 1,\dots, n-1\right\}\right)\setminus \{ 0, n\}  = \left\{2 d m :\, m=1 ,\dots,
\frac{n-d}{2d}\right\}\,,
\end{equation}
and for any
$m\in\left\{1,\dots,\frac{n-d}{2d}\right\} $,
\begin{equation}\label{identity2even}
\card \{ k \in\{1,\ldots,n-1\}:\varphi_j^{(n)}(k)=2 d m\}=2d\,.
\end{equation}}
\end{enumerate}

We prove this lemma establishing a number of intermediate auxiliary results.
\begin{Prop}
\label{Proposition:propsPhi}
Let $n, j \in \N$. Function $\varphi_j^{(n)}:\, \R \to [0,n]$ satisfies:
\begin{enumerate}
\item[(i)] $\varphi_j^{(n)}$ is even: $\varphi_j^{(n)}(-x)=\varphi_j^{(n)}(x)$, $x \in \R$.
\item[(ii)] Symmetry:
\begin{equation}\label{symmetryPhi}
    \varphi_j^{(n)}\left( x+\frac{n}{j}\right) = n- \varphi_j^{(n)}\left( x\right)\,, \quad x \in \R\,.
\end{equation}
In particular, $\varphi_j^{(n)}$ is periodic with period $2n/j$.
\item[(iii)] For every $x\in \R$, either
\begin{equation}\label{mod}
    \frac{\varphi_j^{(n)}(x)-jx }{2n}\in \Z \quad \text{or} \quad \frac{ \varphi_j^{(n)}(x)+jx }{2n}\in \Z\,.
\end{equation}
\item[(iv)] For $k\in \Z$,
\begin{equation}\label{2bis}
\varphi_j^{(n)}(k )=\left|\DD (j k, 2n)\right|\,, \quad k\in \Z\,.
\end{equation}
 
\end{enumerate}
\end{Prop}
 
\begin{proof}
Property \emph{(i)} is obvious from construction. By \eqref{alt1}, if $x \in \left[ \dfrac{2 k-1 }{j}\,
n,\dfrac{2 k+1}{j}\, n\right]$ for $k\in \Z$, then
$$
\varphi_j^{(n)}(x ) =   \begin{cases} 2kn-jx , &  \text{if }
x\in\left[ \dfrac{2k-1}{ j}\, n ,\dfrac{2k}{ j}\, n\right],  \\[4mm]
jx- 2kn  , &  \text{if } x\in\left[ \dfrac{2k}{ j}\, n ,\dfrac{2 k +1}{ j }\, n\right]\,.
\end{cases}
$$
In particular, if $x \in \left[ \dfrac{2 k-1 }{j}\,
n,\dfrac{2 k }{j}\, n\right]$, then $x +n/j \in \left[ \dfrac{2 k  }{j}\,
n,\dfrac{2 k +1}{j}\, n\right]$.  Thus,
$$
\varphi_j^{(n)}(x+n/j ) = j\left( x+\frac{n}{j}\right) - 2kn\,, \quad \varphi_j^{(n)}(x  )= 2kn-jx\,,
$$
and \eqref{symmetryPhi} follows. The case $x \in \left[ \dfrac{2 k  }{j}\,
n,\dfrac{2 k +1 }{j}\, n\right]$ is analyzed analogously.

Identity \eqref{mod} is a straightforward consequence of the explicit formula \eqref{alt1}.

Assume that $k\in \Z$ and $\DD (j k, 2n)=r$; it means that there exists $t\in \Z$ such that
$ r =j k - 2 n t \in [-n, n)  $. Thus, $|\DD (j k, 2n)|= |jk- 2 n t  |$. Furthermore, from the inequalities $-n \leq r < r$ it follows that
$$
\frac{ 2t-1 }{j}\, n \leq k < \frac{ 2t-1 }{j}\, n\,.
$$
Comparing it with the definition of $\varphi_j^{(n)} $ in \eqref{alt1}, we establish \eqref{2bis}.
\end{proof}

Although function $\varphi_j^{(n)}$ is well defined on whole $\R$, we will be mainly interested in its values on the interval $[0,n]$.
In particular, we need to study how $\varphi_j^{(n)}$ acts on integers $1, 2, \dots n-1$:
\begin{Prop}\label{lemmaImportantBis}
Let $n, j\in \N$ with $\text{GCD}(j,n)=1$.

\begin{enumerate}
\item[(i)] For any $j\in \N$, $\varphi_j^{(1)}(\Z)\subset \{0,1 \}$.

\item[(ii)] If $n>1$ and $j$ is odd, then
$$
 \varphi_j^{(n)}  :\, \left\{ 1,\dots, n-1\right\}  \to \left\{ 1,\dots, n-1\right\}
$$
is a bijection. In other words, $\varphi_j^{(n)}$ acts as a permutation
on the set $\{1,\dots,n-1\}$.

\item[(iii)] If $n>1$ and $j$ is even (and thus $n$ is odd), then
$$
 \varphi_j^{(n)} :\, \left\{ 1,\dots, \frac{n-1}{2}\right\}  \to    \left\{ 2m:\, m= 1,\dots, \frac{ n-1}{2}\right\}
$$
and
$$
\varphi_j^{(n)} :\, \left\{ \frac{n+1}{2},\dots, n-1 \right\}  \to  \left\{ 2m:\, m= 1,\dots, \frac{ n-1}{2}\right\}
$$
are bijections.
\end{enumerate}
\end{Prop}
\begin{proof}
By construction, $\varphi_j^{(n)}(k)\in \{ 0, n\}$ if and only if $k= m n/j$, with $m \in \Z $. Since $\text{GCD}(j,n)=1$, $m n/j \in \Z$ only if $m$ is a multiple of $j$. Hence, $\varphi_j^{(n)}(k)\in \{ 0, n\}$ if and only if $k \in \{m n:\, m\in \Z \}$. With $n=1$ this yields \emph{(i)}. Furthermore, if $n>1$, then
$$\varphi_j^{(n)}\left( \{   1, \dots, n-1 \} \right) \subset \{   1, \dots, n -1 \}\,. $$
Hence, in order to prove \emph{(ii)} it is sufficient to show that $\varphi_j^{(n)}$ is injective on $\{   1, \dots, n-1 \} $. Indeed, by \eqref{2bis}, $\varphi_j^{(n)}(x)=r$ if $\left|\DD (j x, 2n)\right|=r$, that is, if there exists $u\in \Z$ such that
$$
j x \pm r = 2n u\,,
$$
with an appropriate choice of the sign. Thus, if $\varphi_j^{(n)}(x)=\varphi_j^{(n)}(y)$, then there exists $u\in \Z$ such that
$$
j (x \pm y) = 2n u\,,
$$
again with an appropriate choice of the sign. However, since $j$ is odd and $\text{GCD}(j,n)=1$, we conclude that $\text{GCD}(j,2n)=1$. This means that $x\pm y$ must be divisible by $2n$. But $|x \pm y|<2n$, so this identity is possible only if $x=y$.

If $j$ is even, then a similar analysis shows that $k \to \left|\DD (j k/2, n)\right|$ is injective both on
$$
\left\{ 1,\dots, \frac{n-1}{2}\right\} \quad \text{and} \quad \left\{ \frac{n+1}{2},\dots, n-1 \right\}\,.
$$
It remains to use that by \eqref{2bis}, $\varphi_j^{(n)}(k )=2 \left|\DD (j k/2, n)\right|$, $ k\in \Z$. This establishes \emph{(iii)}.
\end{proof}

\begin{Prop}
\label{lemma:divisibility}
Let $m , n \in \N$, and $\text{GCD}(m,n)=d$. Then 
$$
\text{GCD}(m,2n)=d \quad \Leftrightarrow \quad m/d \text{ is odd,}
$$
and
$$
\text{GCD}(m,2n)=2d \quad \Leftrightarrow \quad m/d \text{ is even, and } n/d \text{ is odd.}
$$
\end{Prop}
\begin{proof}
Assume that $d=1$. It is obvious that $\text{GCD}(m,2n)=1$ only if $m$ is odd,
and viceversa, if $m$ is odd and $\text{GCD}\left(m,  n \right)=1$, then necessarily $\text{GCD}\left(m,  2n \right)=1$.

Analogously, if $\text{GCD}(m,2n)=2$, it means that $m$ is even, and since $m$ and $n$ are coprime, $n$ must be odd. The reciprocal is also trivially true: if $m$ is even and $\text{GCD}(m, n)=1$, then also $\text{GCD}(m/2,n)=1$, and
$$
 \text{GCD}(m,2n)=2\cdot \text{GCD}(m/2,n)=2\,.
$$
The general case is reduced to $d=1$ by observing that
$$
\text{GCD}(m,n)=d \quad \Leftrightarrow \quad \text{GCD}\left(\frac{m}{d},  \frac{n}{d}\right)=1 \,.
$$
\end{proof}

Now we are ready to prove the main result of this section.

\begin{proof}[Proof of the Main Lemma]

Formula \eqref{hiddenSymbis} is a straightforward consequence of \eqref{2bis}. Since
$$
\varphi_j^{(n)}\left( x+\frac{n}{d}\right) = \varphi_j^{(n)}\left( x+\frac{n}{j}\, \frac{j}{d}\right)\,,
$$
formula \eqref{symmetryPhiMore} follows from \eqref{symmetryPhi}.

Since
$$
\text{GCD}(j,n)=d \quad \Rightarrow \quad \text{GCD}\left(\frac{j}{d},\frac{n}{d}\right)=1\,,
$$
we can apply Proposition \ref{lemmaImportantBis} to function $\varphi_{j/d}^{(n/d)} $. In fact, statement (ii) is a straightforward consequence of \emph{(i)} of Proposition \ref{lemmaImportantBis} and formula \eqref{hiddenSymbis}.

Assume that $\text{GCD}(j,n)=\text{GCD}(j,2n)=d<n$; by Proposition \ref{lemma:divisibility}, $j/d$ is odd. By \emph{(ii)} of Proposition \ref{lemmaImportantBis},
$$
 \varphi_{j/d}^{(n/d)}  :\, \left\{ 1,\dots, n/d-1\right\}  \to \left\{ 1,\dots, n/d-1\right\}
$$
is a bijection, and by formula \eqref{hiddenSymbis}, this is valid also for
$$
\varphi_{j}^{(n)}  :\, \left\{ 1,\dots, n/d-1\right\}  \to \left\{dm:\, m= 1,\dots, n/d-1\right\}\,.
$$
Furthermore, any $k\in \{ 1, \dots, n-1 \}$ can be represented as
$$
k=r+m\frac{n}{d}\,, \quad m\in \{0,1, \dots, d-1 \}\,, \quad r \in \{0, 1, \dots, n/d-1 \}\,.
$$
By \eqref{symmetryPhiMore},
$$
\varphi_j^{(n)}\left( k \right) =\varphi_j^{(n)}\left( r+m\, \frac{n}{d}\right) = \begin{cases}
 n -   \varphi_j^{(n)}\left(r \right)\,, & \text{ if $m$ is odd,} \\
 \varphi_j^{(n)}\left( r \right)\,, & \text{ if $m$ is even.}
\end{cases}
$$
In consequence, for every $m\in \{0,1, \dots, d-1 \}$,
$$
\varphi_{j}^{(n)}  :\, \left\{ 1 + m\, \frac{n}{d}, 2+  m\, \frac{n}{d}, \dots,  (m+1)\, \frac{n}{d}-1 \right\}  \to \left\{dm:\, m= 1,\dots, n/d-1\right\}
$$
is a bijection. This proves \eqref{identity1bis}--\eqref{identity2bis}.

On the other hand, if  $\text{GCD}(j,2n)=2d<2n$, then by Proposition \ref{lemma:divisibility},  $j/d$ is even and $n/d$ is odd, and by \emph{(iii)} of Proposition \ref{lemmaImportantBis},
$$
 \varphi_{j/d}^{(n/d)} :\, \left\{ 1,\dots, \frac{n/d-1}{2}\right\}  \to \left\{ 2m:\, m=1,\dots, \frac{n/d-1 }{2}\right\}
$$
and
$$
\varphi_{j/d}^{(n/d)} :\, \left\{ \frac{n/d+1}{2},\dots, n/d-1 \right\}  \to \left\{ 2m:\, m=1,\dots, \frac{n/d-1 }{2}\right\}
$$
are bijections, so that by formula \eqref{hiddenSymbis}, this is valid also for
$$
 \varphi_{j }^{(n )} :\, \left\{ 1,\dots, \frac{n/d-1}{2}\right\}  \to \left\{ 2dm:\, m=1,\dots, \frac{n/d-1 }{2}\right\}
$$
and
$$
\varphi_{j }^{(n )} :\, \left\{ \frac{n/d+1}{2},\dots, n/d-1 \right\}  \to \left\{ 2dm:\, m=1,\dots, \frac{n/d-1 }{2}\right\}\,.
$$
Again, if $k\in \{ 1, \dots, n-1 \}$, and
$$
k=r+m\frac{n}{d}\,, \quad m\in \{0,1, \dots, d-1 \}\,, \quad r \in \{0, 1, \dots, n/d-1 \}\,,
$$
we have by \eqref{symmetryPhiMore},
$$
\varphi_j^{(n)}\left( k \right) =\varphi_j^{(n)}\left( r+m\, \frac{n}{d}\right) = \varphi_j^{(n)}\left( r\right)\,.
$$
In consequence, for every $m\in \{0,1, \dots, d-1 \}$,
$$
\varphi_{j}^{(n)}  :\, \left\{ 1 + m\, \frac{n}{d}, 2+  m\, \frac{n}{d}, \dots, \frac{n/d-1}{2} + m\, \frac{n}{d}  \right\}  \to \left\{ 2dm:\, m=1,\dots, \frac{n/d-1 }{2}\right\}
$$
and
$$
\varphi_{j}^{(n)}  :\, \left\{ \frac{n/d+1}{2} + m\, \frac{n}{d},   \dots, \frac{n}{d}-1 + m\, \frac{n}{d}  \right\}  \to \left\{ 2dm:\, m=1,\dots, \frac{n/d-1 }{2}\right\}
$$
are bijections. This proves \eqref{identity1even}--\eqref{identity2even}.
\end{proof}

\section{Discrete entropy for Chebyshev polynomials of the first kind} \label{subsection:1kind}

From the explicit formulas \eqref{explicitFirstKind} for $p_n$ is easy to compute that in this case
$$
\ell_n^{-1}(\lambda )=n-\frac{1}{2}+\frac{1}{2}\,\frac{\sin(2n-1)
\theta}{\sin \theta}\,,
$$
and for $\S_n(\lambda )$ defined in \eqref{def_entropy_2} we have
\begin{align*}
 \S_{n }(\lambda ) &= - \log\left(\ell_n(\lambda
)\right) - \ell_n(\lambda  ) \sum_{i=0}^{n-1} \,p_{i }^{2}(\lambda )
\log \left( p_{i }^{2}(\lambda  ) \right) \\
& =- \log\left(\ell_n(\lambda )\right) -
\ell_n(\lambda  ) \sum_{i=1}^{n-1} \,p_{i }^{2}(\lambda ) \log
\left( 2 \cos^{2}(i \theta)(\lambda  ) \right) \\
 &= - \log\left(\ell_n(\lambda )\right) -
\ell_n(\lambda  ) \log  ( 2)\,  \sum_{i=1}^{n-1} \,p_{i
}^{2}(\lambda )  - \ell_n(\lambda ) \sum_{i=1}^{n-1} \,p_{i
}^{2}(\lambda ) \log \left(   \cos^{2}(i \theta)  \right)
\\ =& - \log\left(\ell_n(\lambda )\right) -
\ell_n(\lambda  ) \log  ( 2)\, \left( \ell_n^{-1}(\lambda )-1\right)
- 2 \ell_n(\lambda ) \sum_{i=1}^{n-1} \cos^{2}(i \theta)  \, \log
\left( \cos^{2}(i \theta)  \right)\,,
\end{align*}
so that
\begin{equation}\label{entropyFor_Cheb1}
\S_{n }(\lambda )= - \log\left(\ell_n(\lambda )\right) +   \log  (
2)\, \left( \ell_n (\lambda )-1\right) - 2 \ell_n(\lambda )
\sum_{i=1}^{n-1} \cos^{2}(i \theta) \, \log \left( \cos^{2}(i
\theta)  \right)\,.
\end{equation}

Since for Chebyshev polynomials of the first kind and degree $n$ the
zeros are
\begin{equation*}\label{zeros_of_cheb}
\lambda _j^{(n)}=\cos \left(\frac{(2j-1)\pi }{2n}\right)\,,
\quad j=1, \dots, n\,,
\end{equation*}
we see that $\ell_n(\lambda _j)=1/n$ for $ j=1, \dots, n$. In
particular, by \eqref{entropyFor_Cheb1}, in this case
\begin{equation}\label{entropyFor_Cheb2}
\S_{n,j}=\S_n\left(\lambda _j^{(n)}\right)=\log \left( \frac{n}{2}\right)+\frac{\log
2}{n} - \frac{2}{n}\, \widehat S_{n,j}\,,
\end{equation}
where
\begin{equation}\label{entropyFor_Cheb3}
\widehat S_{n,j}   \isdef 
\sum_{i=1}^{n-1} \cos^{2} \left(  \frac{(2j-1)\pi }{2n}\, i \right) \,
\log \left( \cos^{2}\left( \frac{(2j-1)\pi }{2n}\, i\right)  \right)\,.
\end{equation}

Formulas \eqref{entropyFor_Cheb2}--\eqref{entropyFor_Cheb3} reduce
the computation of $ \S_{n,j}$ to the analysis of the \emph{modified entropy} $\widehat S_{n,j}$. But first we express $\widehat S_{n,j}$ in terms of the auxiliary functions $\varphi_j^{(n)}$ defined by \eqref{alt1}.

\begin{Prop}
\label{lemma:auxiliar} For $j\in \N$,
\begin{equation}\label{permutations_cos}
\left| \cos  \left(\frac{(2j-1)\pi}{2n}\,
x\right)\right| = \left| \cos  \left(\frac{ \pi}{2n}\,
\varphi_{2j-1}^{(n)}(x)\right)\right| \,, \quad x \in \R\,.
\end{equation}
In particular,
\begin{equation}\label{defShatAlt}
\widehat{S}_{n,j} =\sum_{k=1}^{n-1} \cos^2  \left(\frac{ \pi}{2n}\,
\varphi_{2j-1}^{(n)}(k)\right) \log \left( \cos^2  \left(\frac{ \pi}{2n}\,
\varphi_{2j-1}^{(n)}(k)\right)  \right) \,.
\end{equation}
\end{Prop}
\begin{proof}
By \eqref{mod}, given $x\in \R$, either
$$
 \frac{ \varphi_{2j-1}^{(n)}(x) -(2j-1) x}{2n}\in \Z \quad \text{or} \quad \frac{  \varphi_{2j-1}^{(n)}(x)+(2j-1) x }{2n}\in \Z\,.
$$
Hence, there exists $m\in \Z$ such that either
$$
\frac{\pi}{2n}\, \varphi_j^{(n)}(x ) = \frac{(2j-1) \pi}{2n} \,
x + \pi m\,, \quad\text{or} \quad  \frac{\pi}{2n}\,
\varphi_j^{(n)}(x ) = -\frac{(2j-1) \pi}{2n}  \, x + \pi m\,,
$$
and \eqref{permutations_cos} follows.
\end{proof}

Using the arithmetic properties of $\varphi_j^{(n)}$ established above, we can simplify the expression for the modified entropy:
\begin{Prop}
Let $n\in \N$ and $j\in \{ 1, 2, \dots, n\}$. If
$\text{GCD}(2j-1,n)=d $ then
\begin{equation}\label{12}
\widehat{S}_{n,j} = \begin{cases} 0 &
\text{if } j=\dfrac{n+1}{2} \,, \\[3mm]
\displaystyle d
\sum_{k=1}^{(n/d)-1} \cos^2  \left(\frac{\pi  d}{2n}\,
k \right) \log \left( \cos^2  \left(\frac{\pi d}{2n}\,
k \right) \right)\,, & \text{otherwise.}
\end{cases}
\end{equation}
Furthermore,
\begin{equation}\label{symmetryS}
    \widehat{S}_{n,j} = \widehat{S}_{n,n-j+1}\,.
\end{equation}
\end{Prop}
\begin{proof}
Observe that $d\leq n$, and for $j\in \{ 1, 2, \dots, n\}$ we have
$$
d=n \quad \Leftrightarrow \quad j=\dfrac{n+1}{2}\,.
$$
Furthermore, any term with index $k$ in the sum \eqref{defShatAlt}, for which $\varphi_{2j-1}^{(n)}(k )\in \{0, n
\}$, vanishes. Since $\text{GCD}(2j-1,n)= \text{GCD}(2j-1,2n) $, formula \eqref{12} follows in a straightforward way from \eqref{defShatAlt}, (ii)--(iii) of the Main Lemma, and the commutativity of the sum.

Moreover,
\begin{equation}\label{IdentityGCD}
    \text{GCD}(2j-1,n)=\text{GCD}(2(n-j+1)-1,n)\,.
\end{equation}
Indeed, let $\text{GCD}(2j-1,n)=d$, so that $2j-1= d s$, $n=d t$, where $s, t \in \N$  and $\text{GCD}(s,t)=1$. By a well known characterization of coprime integers, there exists integers $x, y $ such that $x s + y t=1$. We have
$ 2(n-j+1)-1= (2t-s) d$, and
$$
- x (2t-s) + (2x + y) t = x s + y t=1 \quad \Rightarrow \quad \text{GCD}(2t-s,t)=1 \,,
$$
so that $\text{GCD}(2(n-j+1)-1,n)=d$, which proves \eqref{IdentityGCD}.
In particular, $\widehat{S}_{n,j}$ and $\widehat{S}_{n,n-j+1}$ share the same value $d$ in \eqref{12}, which proves \eqref{symmetryS}.
\end{proof}

Next we find a series representation for $\widehat S_{n,j}$:
\begin{Prop}
\label{Prop:Scorrected} Let $n\in \N$, $j\in \{ 1, 2, \dots, n\}$,
and $\text{GCD}(2j-1,n)=d$. Then
\begin{equation} \label{seriesForScorrected}
\frac{1}{n}\, \widehat S_{n,j} =    \frac{1}{2 } \, (1-2\log 2)+
\sum_{s=1}^\infty  \zeta ( 2s+1 )\, \left(\frac{d}{2
n}\right)^{2s+1} \,,
\end{equation}
where $\zeta(\cdot)$ is the Riemann zeta function.
\end{Prop}
\begin{remark}
Observe that $d/2n \leq 1/2$, so that the series in the right hand side is convergent.
\end{remark}

\begin{proof} 
In order to find the value of $\widehat S_{n,j}$ we follow the strategy \cite{Aptekarev:95} of computing for $0<2-\epsilon<q< 2+\epsilon$ the $l^q$ norms
\begin{equation}\label{normDef}
N(q;h)\isdef  \sum_{j=0}^{n/d-1}
 \cos^q \left( \frac{\pi d}{2n}\,  j \right)=\sum_{j=0}^{\pi/(2h)-1}
 \cos^q \left( h j \right)\,, \quad h\isdef \frac{\pi d}{2n}\,,
\end{equation}
and considering the partial derivative of $N(q;h)$ with respect to $q$ at $q=2$:
\begin{equation}\label{derivativeofN}\begin{split}
\widehat{S}_{n,j} & =  2d \, \sum_{j=0}^{\pi/(2h)-1} \cos^2\left (h j\right
) \log \cos \left (h j\right) =    2 d  \, \left. \frac{\partial}{\partial q} N(q;h)\right|_{q=2} \,.
\end{split}
\end{equation}
Observe that $N(q;h)$ are related to the Riemann sums of an integral:
\begin{equation}\label{Riemann_sum}
h  N(q;h) \simeq \int_0^{\pi/2} \cos^q (u) du =
\int_0^{\pi/2} \left(\frac{\cos(u)}{ \frac{\pi}{2}-u }\right)^q
\left(\frac{\pi}{2}-u\right)^q du \,.
\end{equation}

Consider $\cos^q (x)$ and $(\cos(x)/(\pi/2-x ))^q$ as analytic functions at $x=0$ and $x=\pi/2$, respectively, whose single valued branches in the corresponding neighborhoods are fixed by
$$
\cos^q (x)\big|_{x=0} = \left(\frac{\cos(x)}{ \frac{\pi}{2}-x }\right)^q \bigg|_{x=\pi/2}= 1.
$$
Denote by
$$
 \cos^q (x)  =  1+\sum_{s=1}^\infty
\alpha_s (q)\, x^s ,\qquad
 \left(\frac{\cos(x)}{ \frac{\pi}{2}-x }\right)^q
  =  1 + \sum_{s=1}^\infty \beta_s (q)\,
\left(\frac{\pi}{2}-x\right)^s
$$
the Taylor expansions of these functions. Observe that $ \alpha_s (q)= \beta_s (q)=0 $ for odd indices $s$.
The Euler-Maclaurin summation formula for integrals with algebraic
singularity at the end points (see \cite{Sidi2004}) yields
\begin{equation}\label{firstApprox}
\begin{split}
h N(q;h) = h \sum_{j=0}^{\pi/(2h)-1}  \cos^q (h j) & = h + \int_0^{\pi/2}
\cos^q (u) du + \\
& + \sum_{s=0}^{\infty} \alpha_s (q)\,  \zeta(-s) \, h^{s+1} +
\sum_{s=0}^\infty \beta_s (q)\, \zeta(-s-q) \, h^{s+q+1}\,;
\end{split}
\end{equation}
at this stage we understand this identity in standard terms of an asymptotic expansion.
Since $ \zeta(0)=-1/2$ and $\zeta(-2j) =0$ for $j\in \N$,
formula in \eqref{firstApprox} reduces to
\begin{equation}\label{secondApprox}
N(q;h) = \frac{1}{2} + \frac{1}{h}\, \int_0^{\pi/2}
\cos^q (u) du +  \sum_{s=0}^\infty \beta_{2s} (q)\, \zeta(-2s-q)
h^{2s+q}\,.
\end{equation}
Paper \cite{Sidi2004} addresses also the case of the Euler-Maclaurin summation formula providing
full asymptotic expansion for integrands with logarithmic
singularities at the end points. In fact, these results from
\cite{Sidi2004} can be obtained by formal differentiation of
\eqref{secondApprox} with respect to $q$. Taking into account \eqref{derivativeofN} we obtain:
\begin{equation}\label{thirdApprox}
\begin{split}
\widehat{S}_{n,j} =& 2 d  \,  \frac{\partial}{\partial q} N(q;h)\big|_{q=2} \\
=& \frac{2 d}{h} \int_0^{\pi/2}\cos^2(u) \log(\cos(u)) du - 2d
 \sum_{s=0}^\infty \beta_{2s} (2) \zeta'(-2(s+1)) h^{2(s+1)}\,.
\end{split}
\end{equation}
But
$$
\int_0^{\pi/2}\cos^2(u) \log(\cos(u)) du = \frac{\pi }{8 }\, (1-\log
4)\,,
$$
and
$$
\frac{\cos^2(z)}{(\frac{\pi}{2}-z)^2} = \sum_{s=0}^{\infty}
\frac{(-1)^s 2^{2s+1}}{(2s+2)!} \left(\frac{\pi}{2}-z\right)^{2s}
\,,
$$
so that
$$
\beta_{2s}(2)   = \frac{(-1)^s 2^{2s+1}}{(2s+2)!}\,.
$$
Recalling the definition of $h$ and gathering these formulas in
\eqref{thirdApprox}, we obtain
\begin{align*} \widehat{S}_{n,j} =&
\frac{n}{2 } \, (1-2\log 2)+ d\, \sum_{s=1}^\infty (-1)^{s} \,
\frac{  \pi^{2s } \zeta'(-2s )}{(2s)!}
\left(\frac{d}{n}\right)^{2s}\,.
\end{align*}
Identity
$$
\zeta'(-2s) = (-1)^s \zeta(2s+1) \frac{(2s)!}{\pi^{2s}}
\frac{1}{2^{2s+1}}
$$
yields now \eqref{seriesForScorrected}. It remains to observe that
the series in the right hand side is convergent, thus this is a bona
fide series expansion of $\widehat{S}_{n,j}$.
\end{proof}

\begin{corollary}
\label{thm1} Let $n\in \N$, $j\in \{ 1, 2, \dots, n\}$, and
$\text{GCD}(2j-1,n)=d$, then
\begin{equation} \label{14}
\frac{2}{n}\, \widehat S_{n,j} =      1-2\log 2  -  \mathcal R
\left(\frac{d}{2 n}\right) \,,
\end{equation}
with $\mathcal R$ defined in \eqref{defR}.
\end{corollary}
\begin{proof}
It is an immediate consequence of \eqref{seriesForScorrected} and \eqref{series_for_R}.
\end{proof}
\begin{remark}
It is easy to check that $\mathcal R(1/2)=1-\log(4)$, so that for
$d=n$ we have $\widehat S_{n,j}=0$ (cf.\ formula \eqref{12}).
\end{remark}

It remains to use \eqref{14}  in \eqref{entropyFor_Cheb2} in order to complete the proof of Theorem \ref{corollary4}.

\section{Entropy of Chebyshev polynomials of the second kind} \label{subsection:2kind}
 
From the explicit formulas \eqref{explicitSecondKind} it follows that the zeros of the Chebyshev polynomials of the second kind of degree $n$ are
$$
\lambda _j^{(n)}=\cos \left( \frac{j \pi }{n+1}\right)\,,
\quad j=1, \dots, n\,,
$$
and
\begin{equation}\label{valueEllSecondKind}
    \ell_n^{-1}(\lambda ) = \sum_{k=0}^{n-1} p_k^2(\lambda )= \frac{1}{2} \frac{n \sin(\theta) -
\cos((n+1)\theta) \sin(n \theta)}{\sin^3(\theta)}, \quad \lambda=\cos\theta\,.
\end{equation}
By \eqref{explicitSecondKind}, with $\lambda =\cos (\theta)$,
\begin{align*}
 \sum_{k=1}^{n} p_{k-1}^2\left( \lambda   \right) \log\left(  p_{k-1}^2\left( \lambda   \right)\right) & = \sum_{k=1}^{n} p_{k-1}^2\left( \lambda   \right) \log \left( \sin^2\left( k \theta   \right)\right) -  \log\left(  \sin^2 (\theta)\right)\,   \sum_{k=1}^{n} p_{k-1}^2\left( \lambda   \right) \\
  & = \sin^{-2} (\theta)\, \sum_{k=1}^{n} \sin^2\left( k \theta   \right) \log \left( \sin^2\left( k \theta   \right)\right) -  \log\left(  \sin^2 (\theta)\right)\,  \ell_n^{-1}(\lambda )\,.
\end{align*}
Introducing the notation
\begin{align} \label{defSHatSecondKind}
\widehat{S}_{n,j} & \isdef  
\sum_{k=1}^{n-1} \sin^2\left (\frac{k j\pi}{n } \right )
\log\left ( \sin^2\left ( \frac{k j\pi}{n }\right)\right)   \,,
\end{align}
and using \eqref{def_entropy_1} we get
\begin{align}
\S_{n,j} & =  - \log(\ell_n(\lambda_j^{(n)})) - \ell_n(\lambda_j^{(n)})
\left( \frac{\widehat{S}_{n+1,j} }{\sin^2(j \pi /(n+1))} -
 \frac{\log \left( \sin^2(j \pi /(n+1))\right)}{\ell_n(\lambda_j)}
\right) \nonumber \\ & = - \log(\ell_n(\lambda_j^{(n)})) + \log \left( \sin^2(j \pi /(n+1))\right) -   \frac{\ell_n(\lambda_j^{(n)})  }{\sin^2(j \pi /(n+1))} \, \widehat{S}_{n+1,j} \,. \label{formulaforSsecondkind}
\end{align}
Furthermore, by \eqref{valueEllSecondKind},
\begin{equation}\label{EllSinJ}
    \ell_n^{-1}(\lambda_j^{(n)})=\frac{n+1}{2\, \sin^2 \left(\frac{j \pi }{n+1} \right)}\,,
\end{equation}
so that
\begin{equation}\label{entropyFor_Cheb2Kind}
\S_{n,j}=\S_n\left(\lambda _j^{(n)}\right)=\log \left( \frac{n+1}{2}\right) - \frac{2}{n+1}\, \widehat S_{n+1,j}\,.
\end{equation}
Again, formula \eqref{entropyFor_Cheb2Kind} reduces
the computation of $ \S_{n,j}$ to the analysis of the \emph{modified entropy} $\widehat S_{n,j}$. We express $\widehat S_{n,j}$ in terms of the auxiliary functions $\varphi_j^{(n)}$ defined by \eqref{alt1}:
\begin{Prop}
\label{lemma:auxiliarSecondKind} For $j\in \N$, 
\begin{equation}\label{permutations_sin}
\left| \sin \left(\frac{j \pi}{ n}\,
x\right)\right| = \left| \sin \left(\frac{  \pi}{2 n}\,
\varphi_{2j}^{(n)}(x) \right)\right| \,, \quad x \in \R\,.
\end{equation}
In particular,
\begin{equation}\label{defShatAltSin}
\widehat{S}_{n,j} =\sum_{k=1}^{n-1} \sin^2 \left(\frac{  \pi}{2 n}\,
\varphi_{2j}^{(n)}(k) \right) \log \left( \sin^2 \left(\frac{  \pi}{2 n}\,
\varphi_{2j}^{(n)}(k) \right) \right) \,.
\end{equation}
\end{Prop}
\begin{proof}
By \eqref{mod}, for every $x\in \R$ either
$$
 \frac{ \varphi_{2j}^{(n)}(x) -2j x}{2n}\in \Z \quad \text{or} \quad \frac{  \varphi_{2j}^{(n)}(x)+2j x }{2n}\in \Z\,.
$$
Hence, given $x\in \R$ there exists $m\in \Z$ such that either
$$
\frac{\pi}{2n}\, \varphi_{2j}^{(n)}(x ) = \frac{j \pi}{ n} \,
x + \pi m\,, \quad\text{or} \quad  \frac{\pi}{2n}\,
\varphi_j^{(n)}(x ) = -\frac{j \pi}{ n}  \, x + \pi m\,,
$$
which implies \eqref{permutations_sin}.
\end{proof}

\begin{Prop}
Let $n\in \N$ and $j\in \{ 1, 2, \dots, n-1\}$. Then
\begin{equation}\label{12even}
\widehat{S}_{n,j} = \begin{cases} 0 &
\text{if } j=n/2\,, \\[3mm]
\displaystyle D \, \sum_{m=1}^{\frac{n }{ D}-1}
\sin^2 \left(\frac{  \pi d}{  n}\,
m \right) \log \left( \sin^2 \left(\frac{  \pi d}{  n}\,
m \right) \right),
& \text{otherwise,}
\end{cases}  
\end{equation}
where $D=\text{GCD}(2j ,n)  $ and $d=\text{GCD}( j ,n) $.

Furthermore,
\begin{equation}\label{symmetrySeven}
    \widehat{S}_{n,j} = \widehat{S}_{n,n-j }\,.
\end{equation}
\end{Prop}
\begin{proof}
Observe that $D\leq n$, and for $j\in \{ 1, 2, \dots, n-1\}$ we have
$$
D=n \quad \Leftrightarrow \quad j=\dfrac{n}{2}\,.
$$
It follows from (ii) of the Main Lemma (Section \ref{subsection:endomorphism}) that in this case $\varphi_{2j}^{(n)}(\Z )\in \{0, n
\}$, so that we obtain formula \eqref{12even} for $j=n/2$.

For $j \neq n/2$ we consider two cases. First, assume that $D=2d$.
Then,
$$
\text{GCD}(2j ,2n)= 2 \cdot  \text{GCD}(j ,n) = 2 d = D  =  \text{GCD}(2j ,n) \,.
$$
Since $D<n$, by \eqref{identity1bis}--\eqref{identity2bis}, $
 \varphi_{2j}^{(n)} \left( \left\{ 1,\dots, n-1\right\}\right)\setminus \{ 0, n\}  = \left\{   D m :\, m=1 ,\dots,
\frac{n }{D}-1\right\}$,
and for any $m\in\left\{1,\dots,\frac{n }{ D}-1\right\} $, $ \card \{ k \in\{1,\ldots,n-1\}:\varphi_{2j}^{(n)}(k)=  D m\}= D $. Hence, \eqref{defShatAltSin} and the commutativity of the sum yield
$$
\widehat{S}_{n,j} = D \, \sum_{m=1}^{\frac{n }{ D}-1}
\sin^2 \left(\frac{  \pi D}{ 2 n}\,
m \right) \log \left( \sin^2 \left(\frac{  \pi D}{ 2 n}\,
m \right) \right),
$$
which proves  \eqref{12even} in this case.

Assume next that  $D=  d$.
Then,
$$
\text{GCD}(2j ,2n)= 2 \cdot  \text{GCD}(j ,n) = 2 d = 2 D = 2 \cdot \text{GCD}(2j ,n) \,.
$$
Since $D<n$, by \eqref{identity1even}--\eqref{identity2even}, $
 \varphi_{2j}^{(n)} \left( \left\{ 1,\dots, n-1\right\}\right)\setminus \{ 0, n\}  = \left\{2 D m :\, m=1 ,\dots,
\frac{n/D-1}{2}\right\}$,
and for any $m\in\left\{1,\dots,\frac{n-D}{2D}\right\} $, $ \card \{ k \in\{1,\ldots,n-1\}:\varphi_{2j}^{(n)}(k)=2D m\}=2D $. Using \eqref{defShatAltSin} and the commutativity of the sum, we conclude that
\begin{equation}\label{WideS1}
    \widehat{S}_{n,j} =2D\, \sum_{m=1}^{\frac{n/D-1}{2}} \sin^2 \left(\frac{  \pi D}{  n}\,
m \right) \log \left( \sin^2 \left(\frac{  \pi D}{  n}\,
m \right) \right)\,.
\end{equation}
Furthermore,
$$
\sin^2 \left(\frac{  \pi D}{  n}\,
m \right) = \sin^2 \left( \pi - \frac{  \pi D}{  n}\,
m \right) = \sin^2 \left(\frac{  \pi D}{  n}\,
\left(\frac{  n}{  D} -m \right) \right) \,,
$$
so that
$$
\widehat{S}_{n,j} =2D\, \sum_{m=1}^{\frac{n/D-1}{2}} \sin^2 \left(\frac{  \pi D}{  n}\,
\left(\frac{  n}{  D} -m \right) \right) \log \left( \sin^2 \left(\frac{  \pi D}{  n}\,
\left(\frac{  n}{  D} -m \right) \right) \right)\,.
$$
But
$$
\left\{  m = 1, \dots, \frac{n/D-1}{2} \right\} \cup \left\{ \frac{n}{D}-m :\, m = 1, \dots, \frac{n/D-1}{2} \right\}=\left\{ 1, \dots, \frac{n }{D}-1 \right\}\,.
$$
Thus,
$$
 \widehat{S}_{n,j} = D\, \sum_{m=1}^{\frac{n}{D}-1} \sin^2 \left(\frac{  \pi D}{  n}\,
m \right) \log \left( \sin^2 \left(\frac{  \pi D}{  n}\,
m \right) \right)\,,
$$
which concludes the proof of \eqref{12even}.

Finally, we have that $j= d s$, $n=d t$, where $s, t \in \N$  and $\text{GCD}(s,t)=1$. Again by a characterization of coprime integers, there exists integers $x, y $ such that $x s + y t=1$. But
$ n-j= (t-s) d$, and
$$
- x (t-s) + (x + y) t = x s + y t=1 \quad \Rightarrow \quad \text{GCD}(t-s,t)=1 \,,
$$
so that $\text{GCD}(n-j,n)=d$. Analogously, $\text{GCD}(2j,n)=D=\text{GCD}(2(n-j),n)$.
Now \eqref{symmetrySeven} is a straightforward consequence of \eqref{12even}.
\end{proof}

Next we find a series representation for $\widehat S_{n,j}$:
\begin{Prop}
\label{Prop:ScorrectedSin} Let $n\in \N$, $j\in \{ 1, 2, \dots, n-1\}$,
and $\text{GCD}( j ,n)=d$. Then
\begin{equation} \label{seriesForScorrectedSin}
\frac{2}{n}\, \widehat S_{n,j} =   1-2\log 2  -  \mathcal R
\left(\frac{d}{  n}\right) \,,
\end{equation}
where $\mathcal R(\cdot)$ has been introduced in \eqref{defR}.
\end{Prop}
\begin{proof} 
Observe first that
\begin{align}\nonumber
    \sum_{m=1}^{n-1} \sin^2 \left(\frac{  \pi  }{2 n}\,
m \right) & \log \left( \sin^2 \left(\frac{  \pi  }{ 2 n}\,
m \right) \right) \\ & = \sum_{m=1}^{n-1} \sin^2 \left(\frac{  \pi  }{2 n}\,
(n-m) \right) \log \left( \sin^2 \left(\frac{  \pi  }{ 2 n}\,
(n-m) \right) \right) \nonumber \\  & = \sum_{m=1}^{n-1} \cos^2 \left(\frac{  \pi  }{2 n}\,
m \right) \log \left( \cos^2 \left(\frac{  \pi  }{ 2 n}\,
m \right) \right)\,. \label{sin-cos}
\end{align}
Let $\text{GCD}( 2j ,n)=D$; assume first that $D=2d$. By \eqref{12even} and  \eqref{sin-cos},
\begin{align*}
\widehat{S}_{n,j} & = D \, \sum_{m=1}^{\frac{n }{ D}-1}
\sin^2 \left(\frac{  \pi d}{  n}\,
m \right) \log \left( \sin^2 \left(\frac{  \pi d}{  n}\,
m \right) \right) \\ & =  D \, \sum_{m=1}^{\frac{n }{ D}-1}
\sin^2 \left(\frac{  \pi D}{ 2 n}\,
m \right) \log \left( \sin^2 \left(\frac{  \pi D}{ 2 n}\,
m \right) \right) \\ & = D \, \sum_{m=1}^{\frac{n }{ D}-1}
\cos^2 \left(\frac{  \pi D}{ 2 n}\,
m \right) \log \left( \cos^2 \left(\frac{  \pi D}{ 2 n}\,
m \right) \right) \\
& = \frac{n}{2}\, \left(     1-2\log 2  -  \mathcal R
\left(\frac{D}{2 n}\right) \right)\,.
\end{align*}
where for the last identity we have used \eqref{12} and \eqref{14}. Since $D=2d$, this proves \eqref{seriesForScorrectedSin} in this case.

The remaining case is analyzed in a similar fashion.

\end{proof}

Using \eqref{seriesForScorrectedSin}  in \eqref{entropyFor_Cheb2Kind} we complete the proof of Theorem \ref{corollary42Kind}.

\section{Further numerical experiments} \label{sec5}

In this section we present some results of numerical evaluation of the entropy $\S_{n,j}$ for several orthogonal polynomials. Computation has been carried out in Fortran 95, by complete diagonalization of the corresponding Jacobi matrix $L_n$ in \eqref{JacobiMatrix}, using the routine STEVD of the LAPACK95 library \cite{Anderson99, Barker01}, which computes all the eigenvalues and eigenvectors of a given matrix by means of a divide and conquer algorithm \cite{Rutter94}.

\begin{figure}[htb]
\hspace{-1cm}
\begin{tabular}{ll}
\vspace{-1.2cm}\mbox{\includegraphics[scale=0.4]{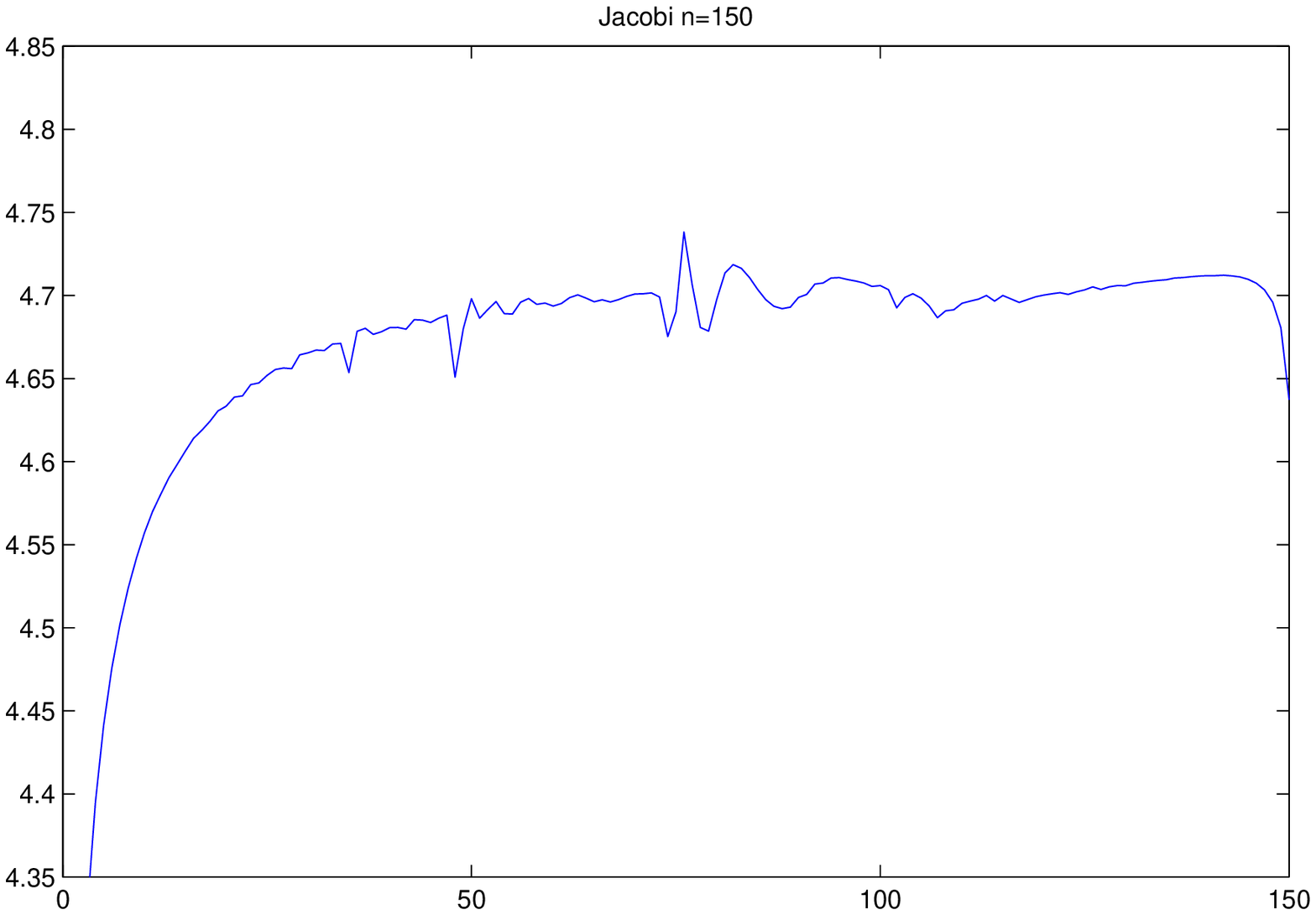}}  & \hspace{-1.2cm}\mbox{\includegraphics[scale=0.4]{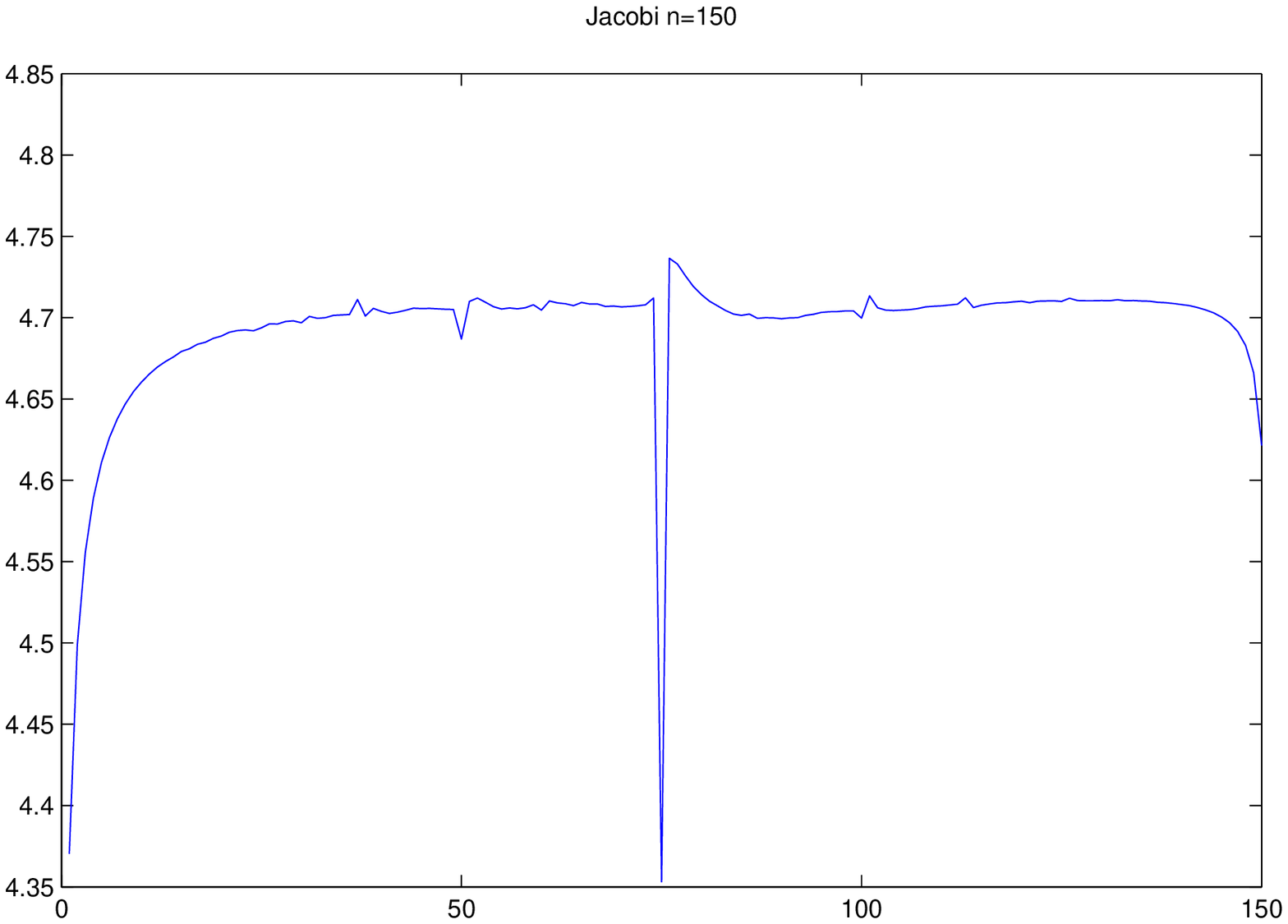}} \\
\vspace{-1cm}\mbox{\includegraphics[scale=0.4]{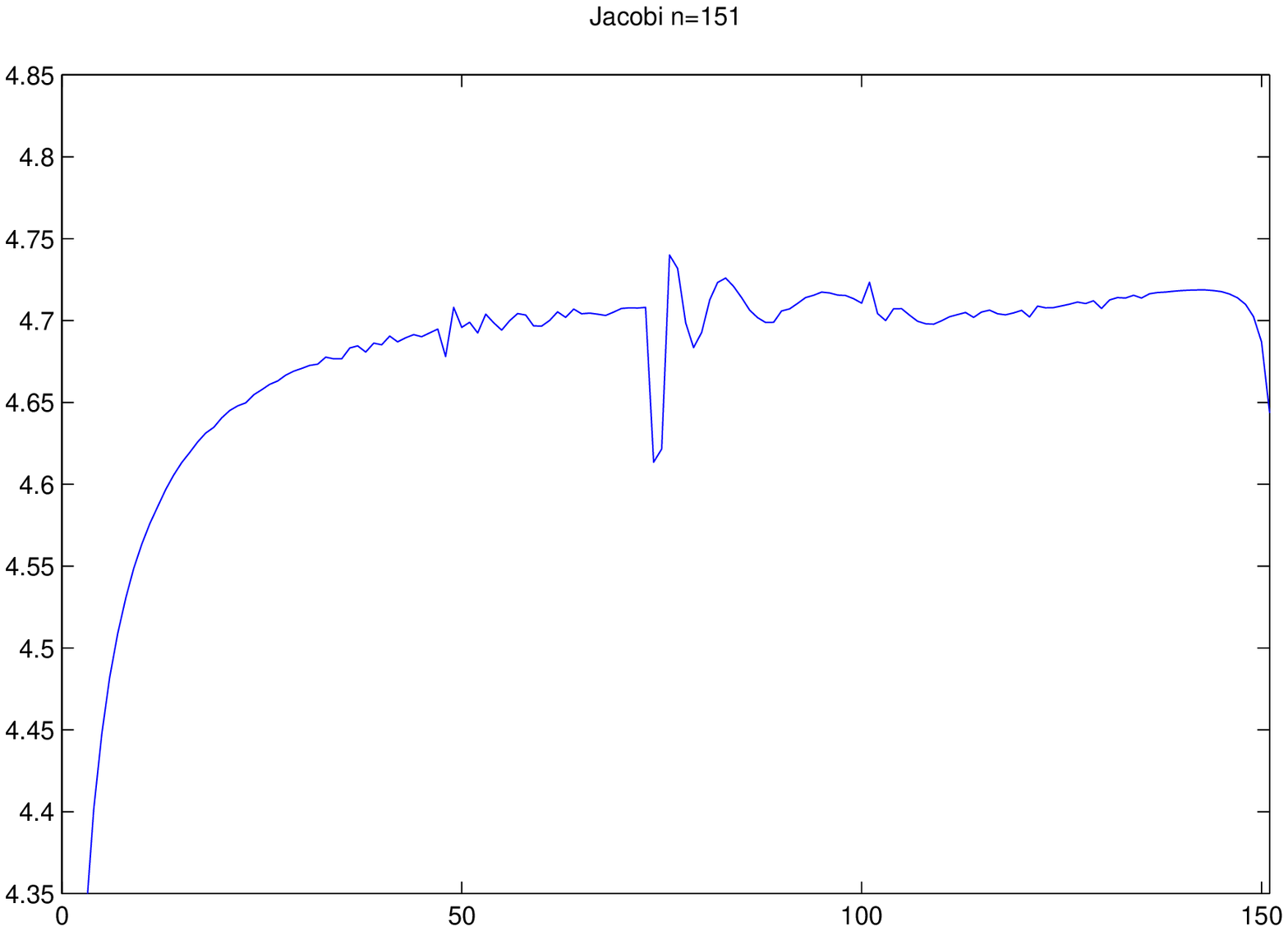}} &
\hspace{-1.2cm}\mbox{\includegraphics[scale=0.4]{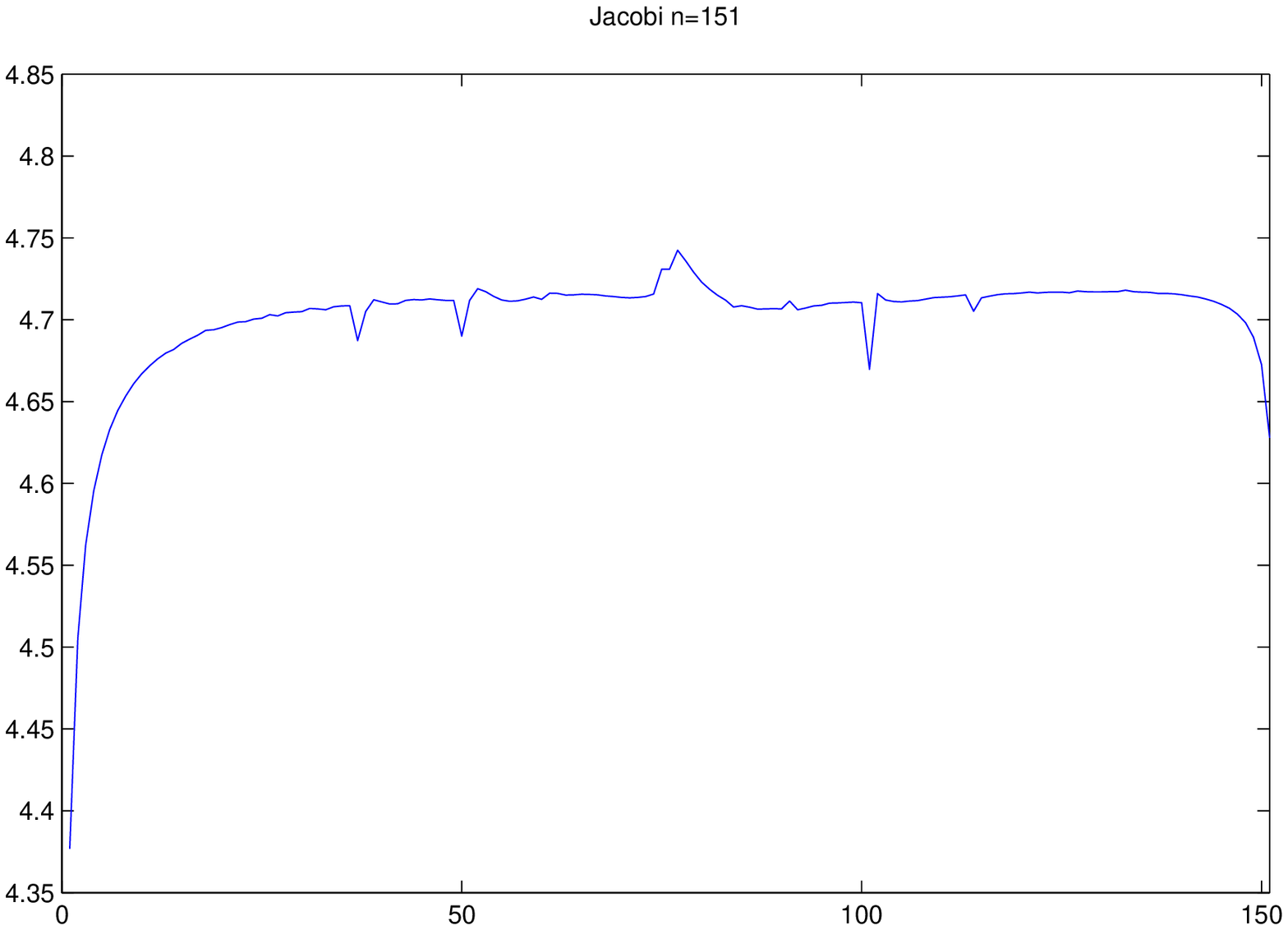}}\\
\vspace{-0.7cm}\mbox{\includegraphics[scale=0.4]{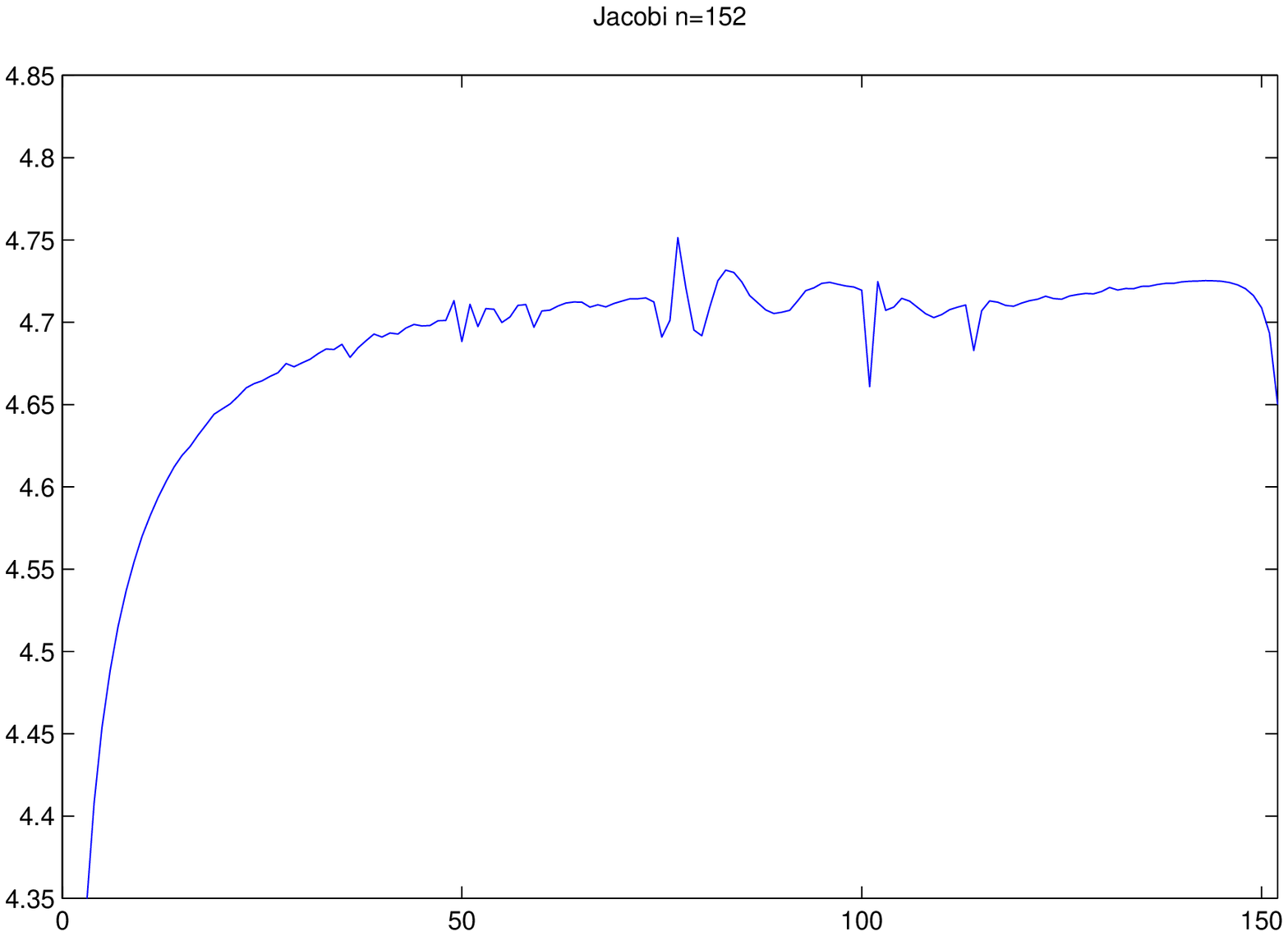}} &
\hspace{-1.2cm}\mbox{\includegraphics[scale=0.4]{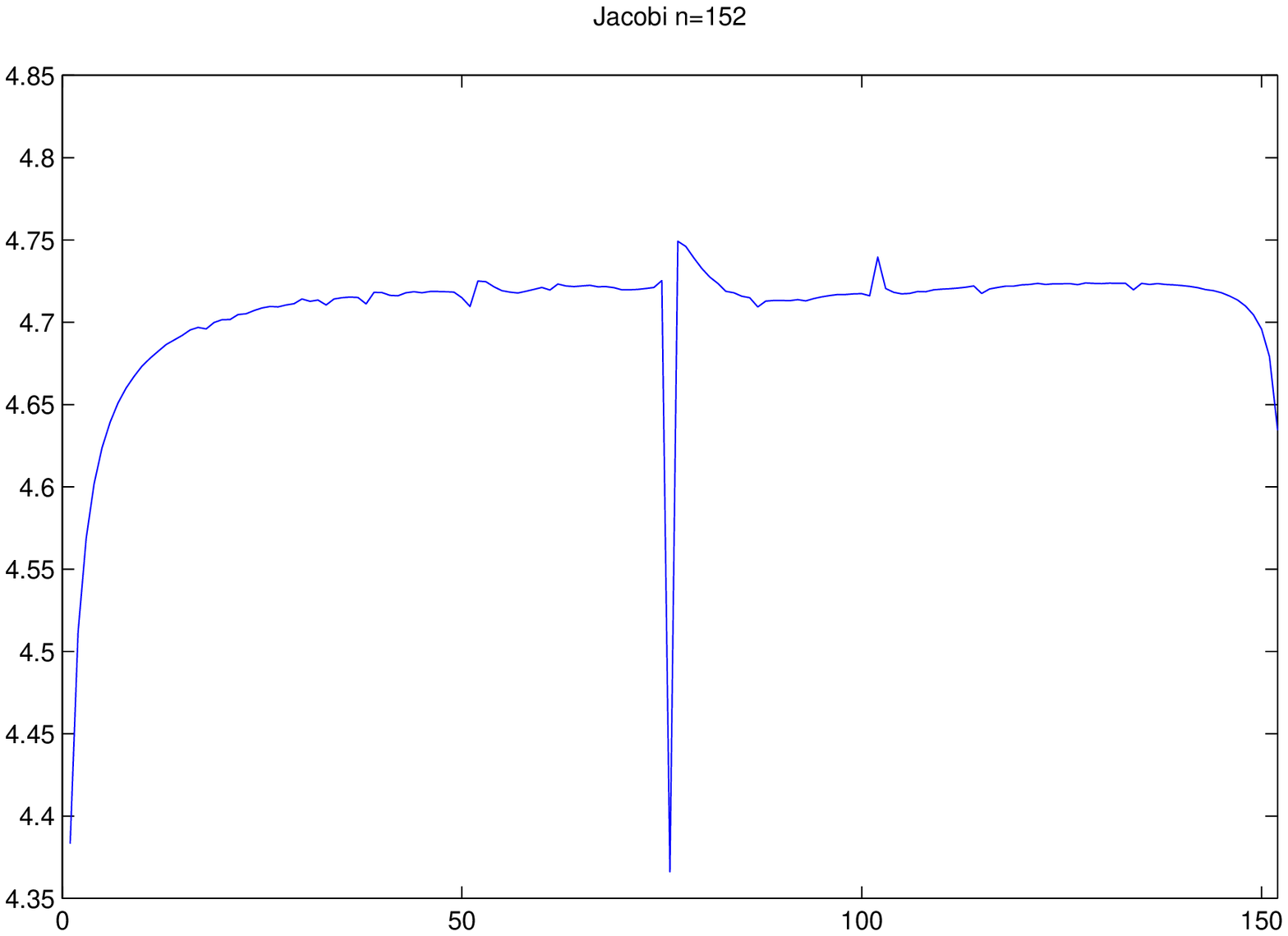}}
\end{tabular}
 \caption{Jacobi polynomials $P_n^{(1.2, 8.9)}$ (left) and $P_n^{(1.2, 3.4)}$ (right): entropy $\S_{n,j}$ for $n=150$, $151$ and $152$.
}\label{fig:Jacobi1}
\end{figure}

As a first illustration we present the entropies $\S_{n, j}$, $n=150, 151, 152$, for two values of the parameters $\alpha $, $\beta $ of Jacobi polynomials $P_n^{(\alpha,\beta)} $, given by the recurrence relation \eqref{recurrence} with
$$
b_i = \frac{2}{2i+\alpha+\beta} \sqrt{\frac{i (i+\alpha) (i+\beta) (i+\alpha+\beta)}{(2i+\alpha+\beta+1)  (2i+\alpha+\beta -1)}}\,, \quad a_i = \frac{\alpha^2-\beta^2}{(2i+\alpha+\beta) (2i+\alpha+\beta-2)}\,.
$$
In Figure \ref{fig:Jacobi1} we can also observe the ``peaks'' explained for the Chebyshev polynomials, but unlike in the latter case, they are pointing both downwards and upwards. Furthermore, the value distribution close the endpoints of the interval clearly differs from the behavior in the bulk.
\begin{figure}[htb]
\hspace{-1cm}
\begin{tabular}{ll}
\vspace{-1.2cm}\mbox{\includegraphics[scale=0.4]{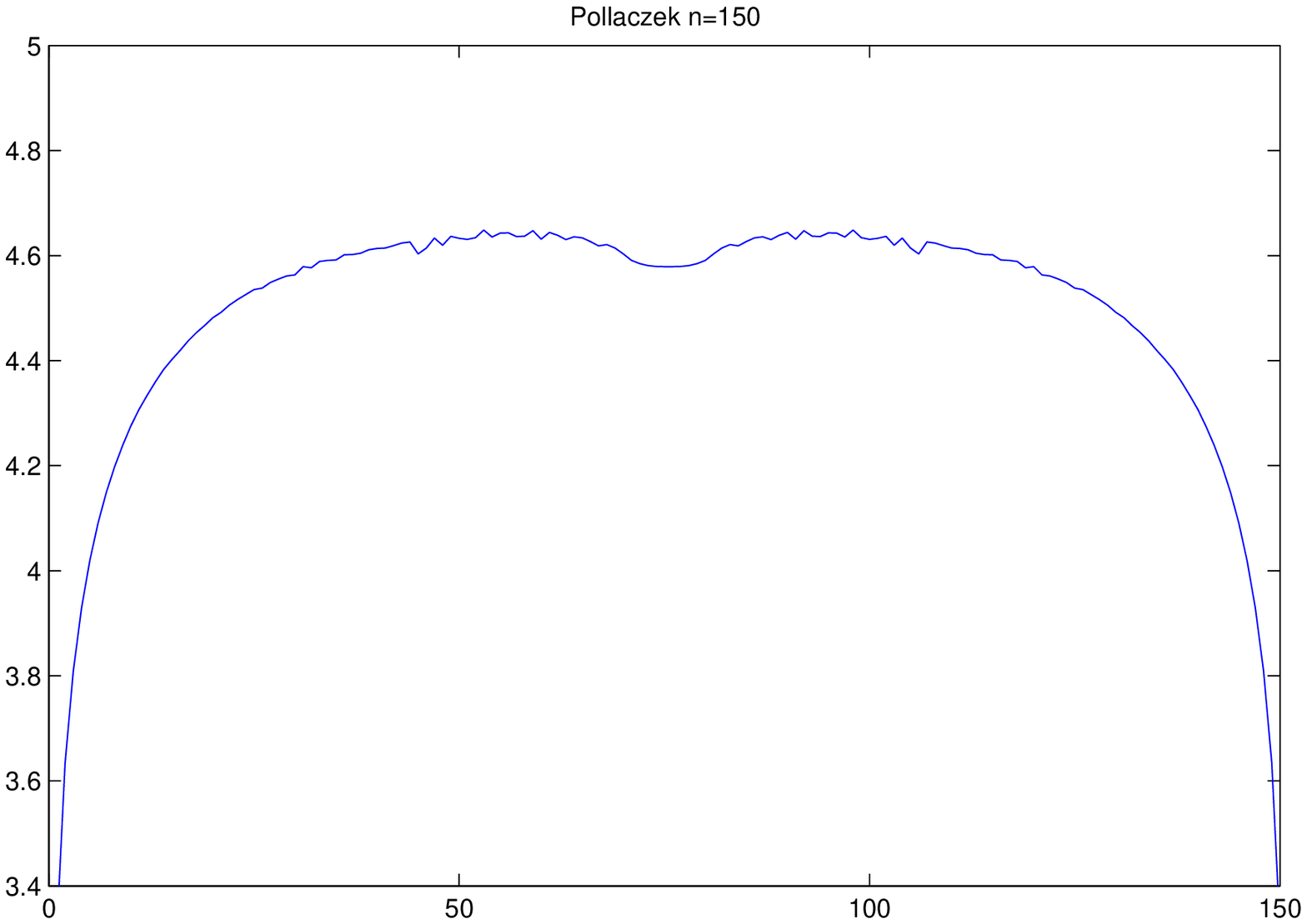}}  & \hspace{-1.2cm}\mbox{\includegraphics[scale=0.4]{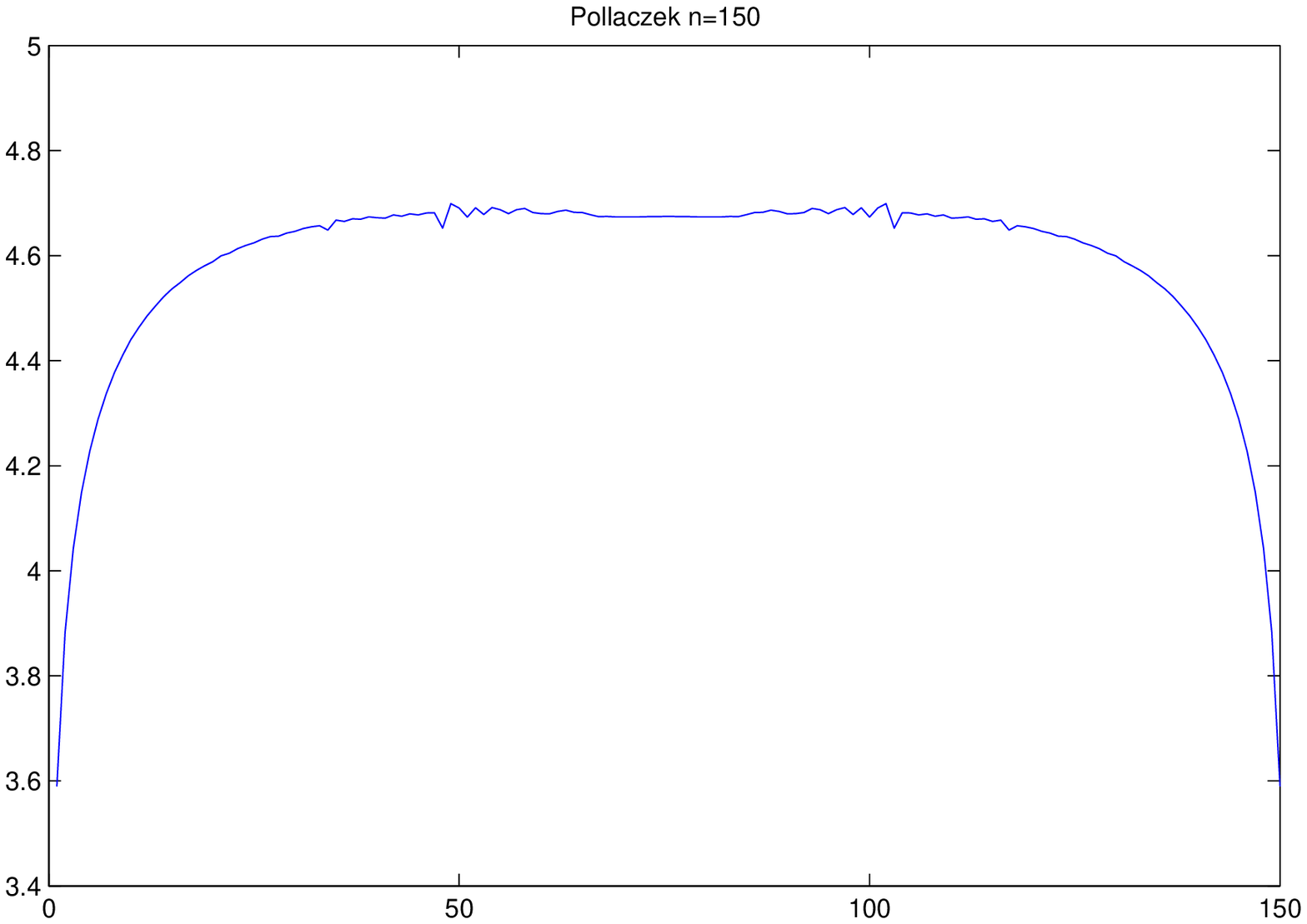}} \\
\vspace{-1cm}\mbox{\includegraphics[scale=0.4]{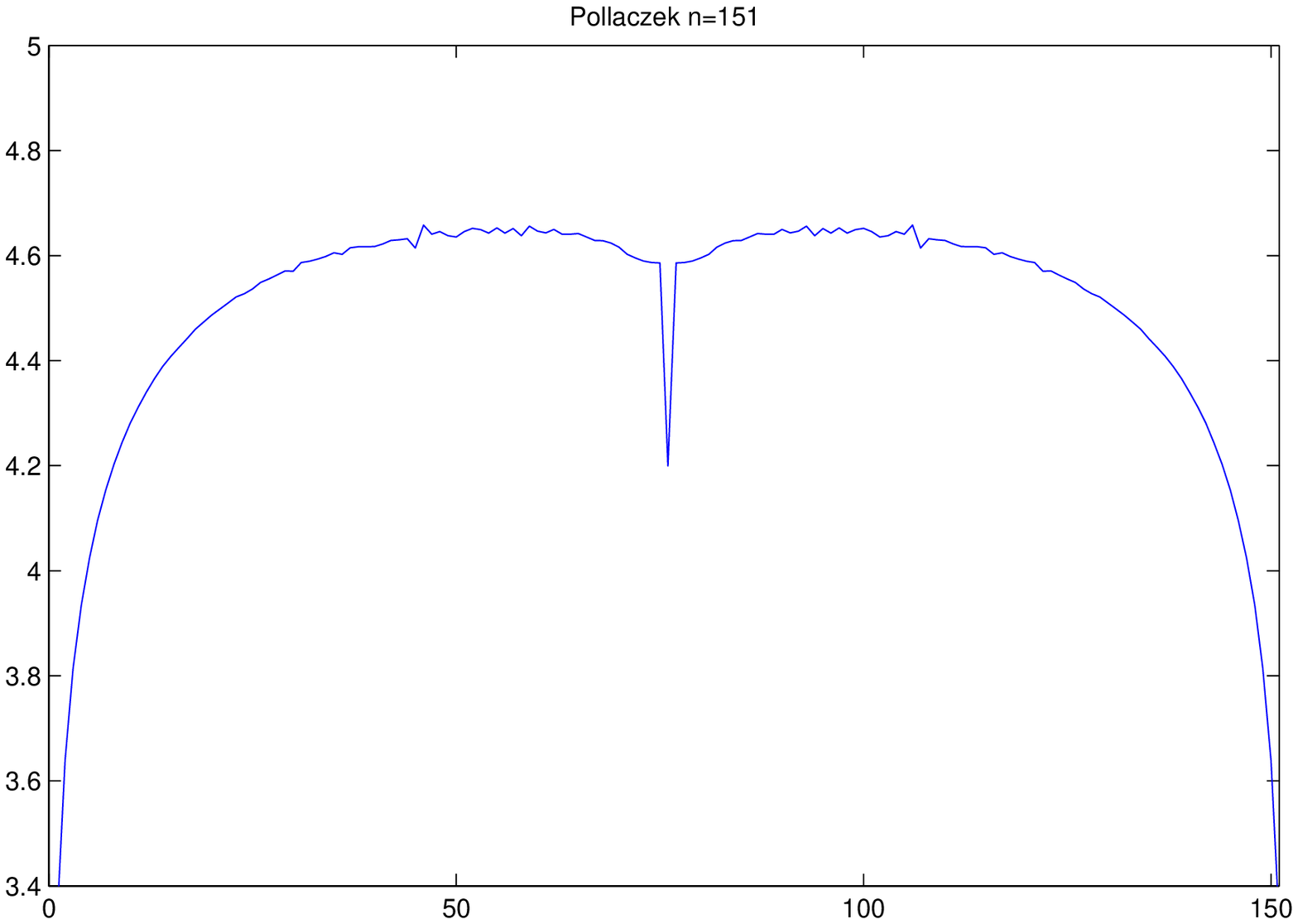}} &
\hspace{-1.2cm}\mbox{\includegraphics[scale=0.4]{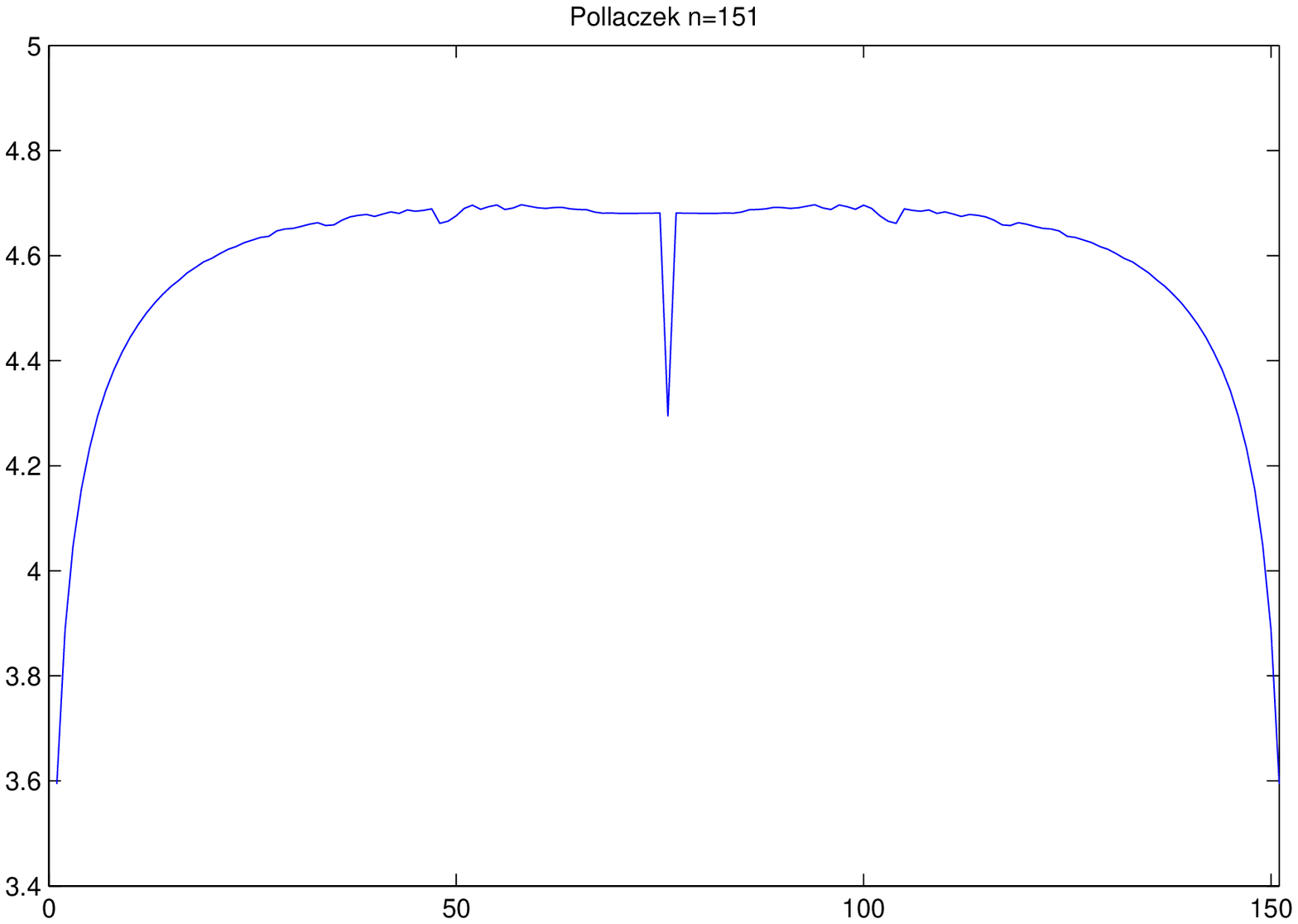}}\\
\vspace{-0.7cm}\mbox{\includegraphics[scale=0.4]{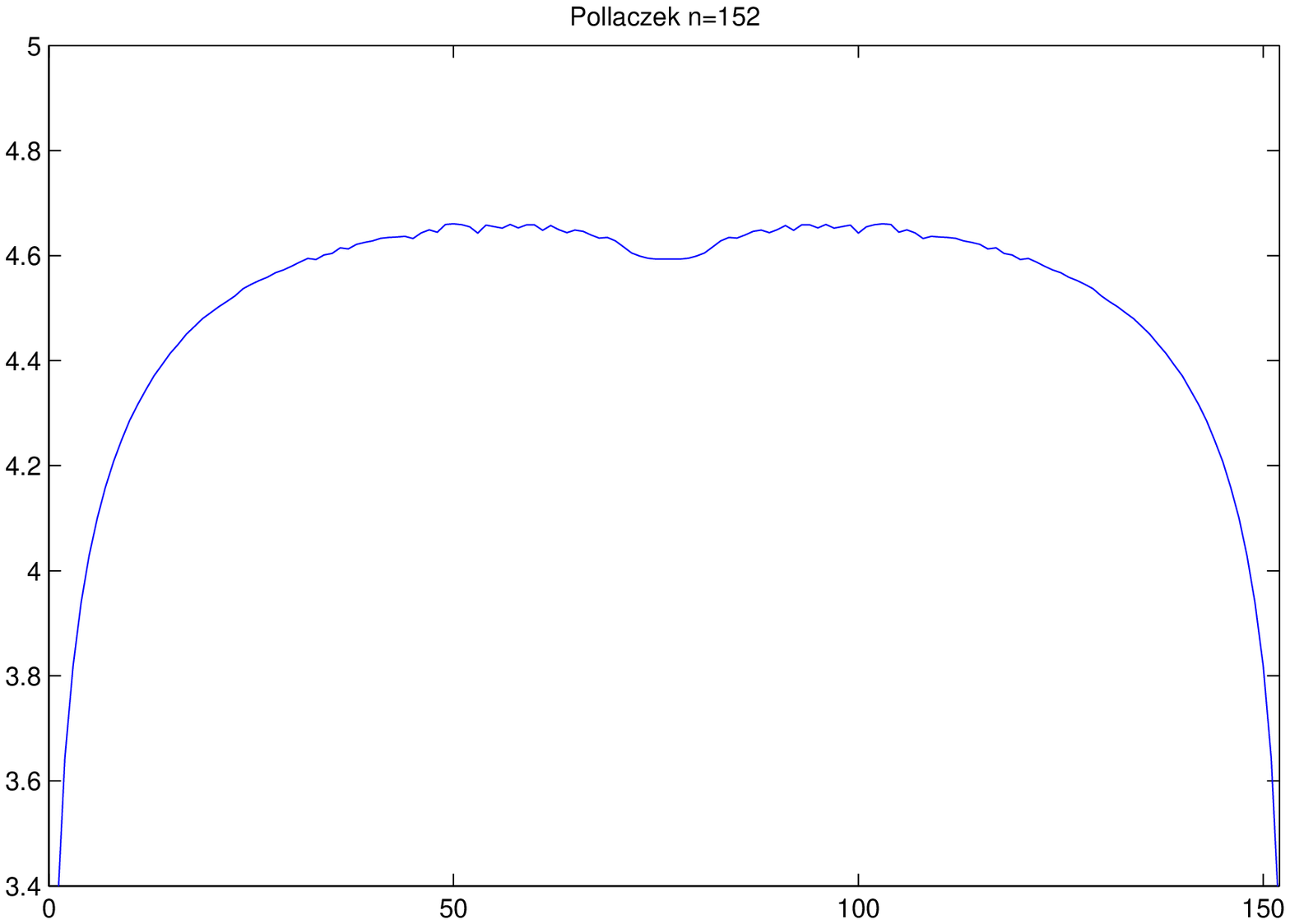}} &
\hspace{-1.2cm}\mbox{\includegraphics[scale=0.4]{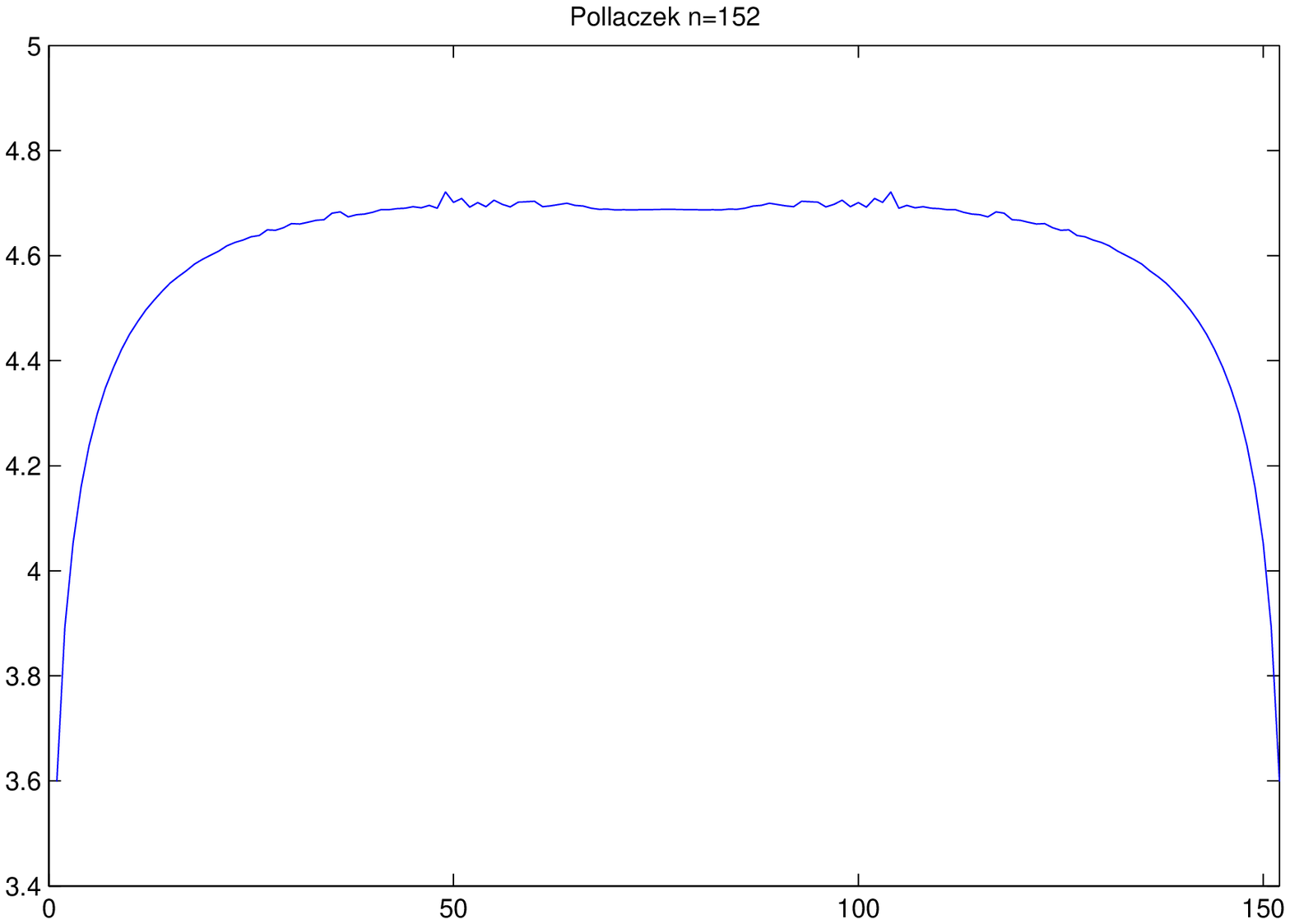}}
\end{tabular}
 \caption{Pollazcek polynomials $p_n^{1.2 }(x;  8.9)$ (left) and $p_n^{1.2  }(x; 3.4)$ (right): entropy $\S_{n,j}$ for $n=150$, $151$ and $152$.
}\label{fig:Pollaczek}
\end{figure}
The feature of the endpoint behavior is even more visible for the symmetric Pollaczek polynomials $p_n^\theta (\cdot ;a)$ (Figure \ref{fig:Pollaczek}), given by the recurrence relation \eqref{recurrence} with
$$
b_i = \frac{1}{2} \sqrt{\frac{i(i+2\theta-1)}{(i+\theta+a)(i+\theta+a-1)}}\,, \quad a_i =0 \,.
$$
For $a =0$ these polynomials reduce to the Jacobi (or more precisely, Gegenbauer) polynomials $P_n^{(\theta -1/2,\theta -1/2)} $. If $a >0$, the orthogonality weight for the Pollaczek polynomials does not satisfy the Szeg\H{o} condition due precisely to its exponentially fast decay at the end points of the interval $[-1, 1]$.

\begin{figure}[htb]
\hspace{-1cm}
\begin{tabular}{ll}
\vspace{-0.9cm}\mbox{\includegraphics[scale=0.4]{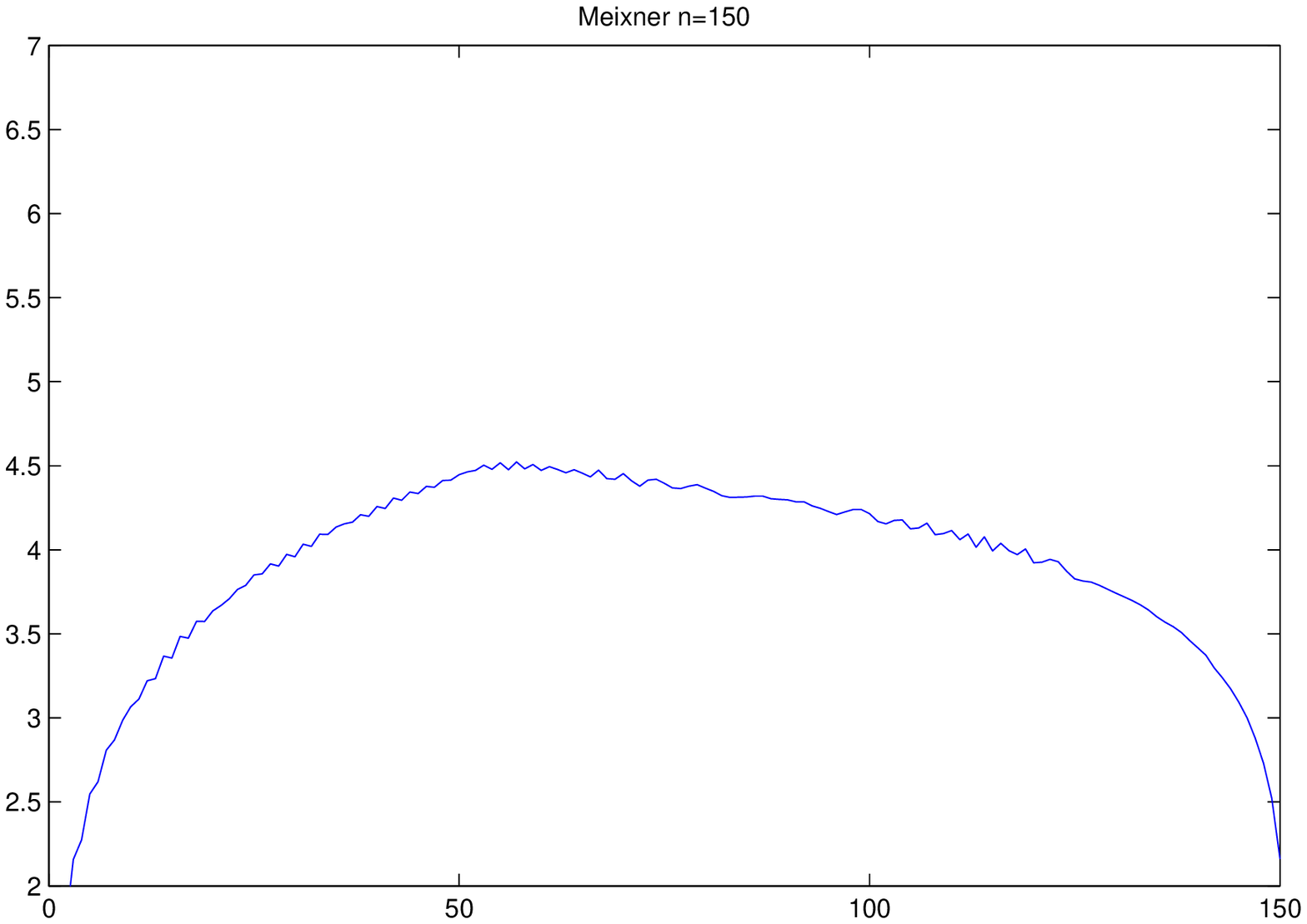}}  & \hspace{-1.2cm}\mbox{\includegraphics[scale=0.4]{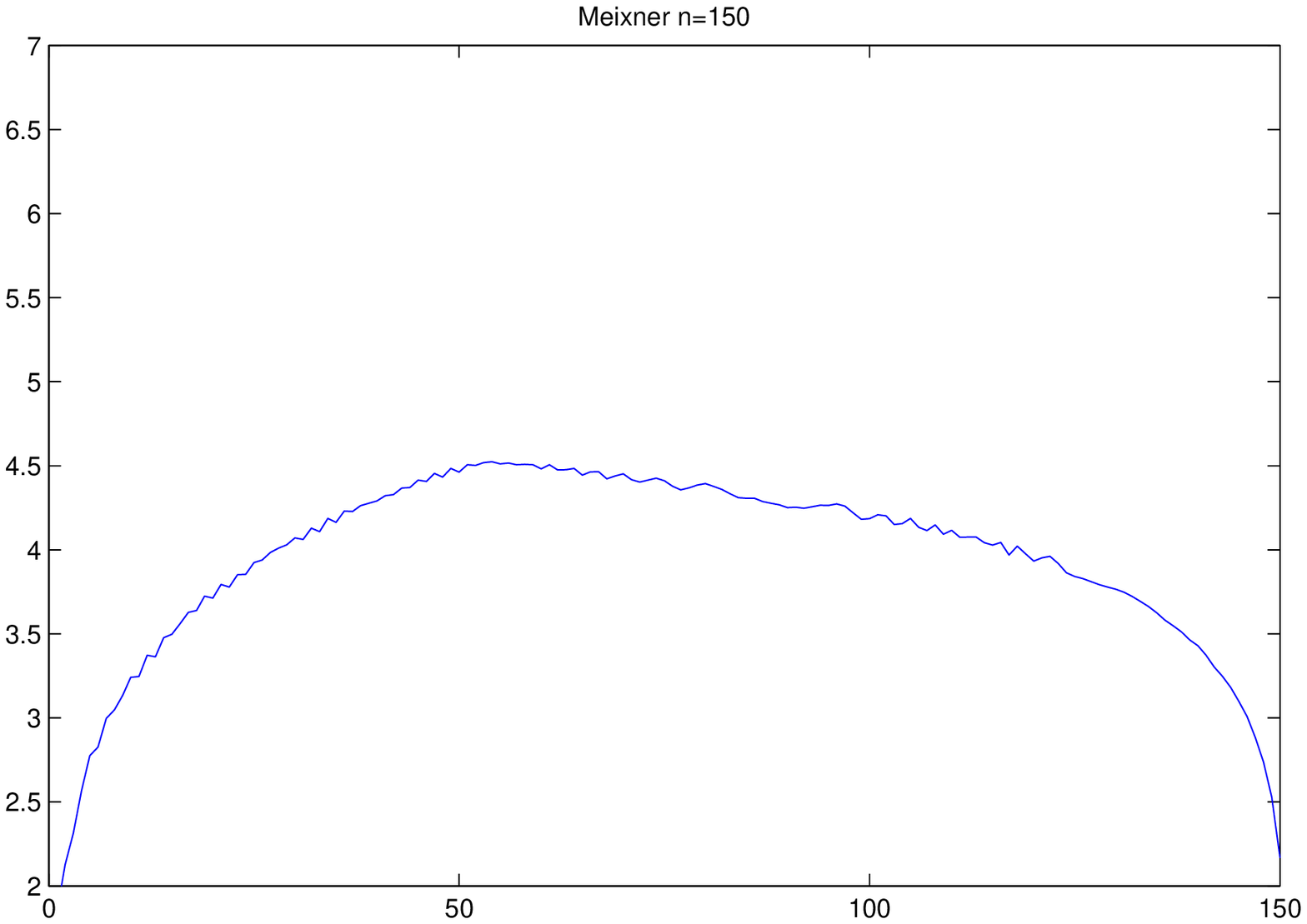}} \\
\vspace{-0.7cm}\mbox{\includegraphics[scale=0.4]{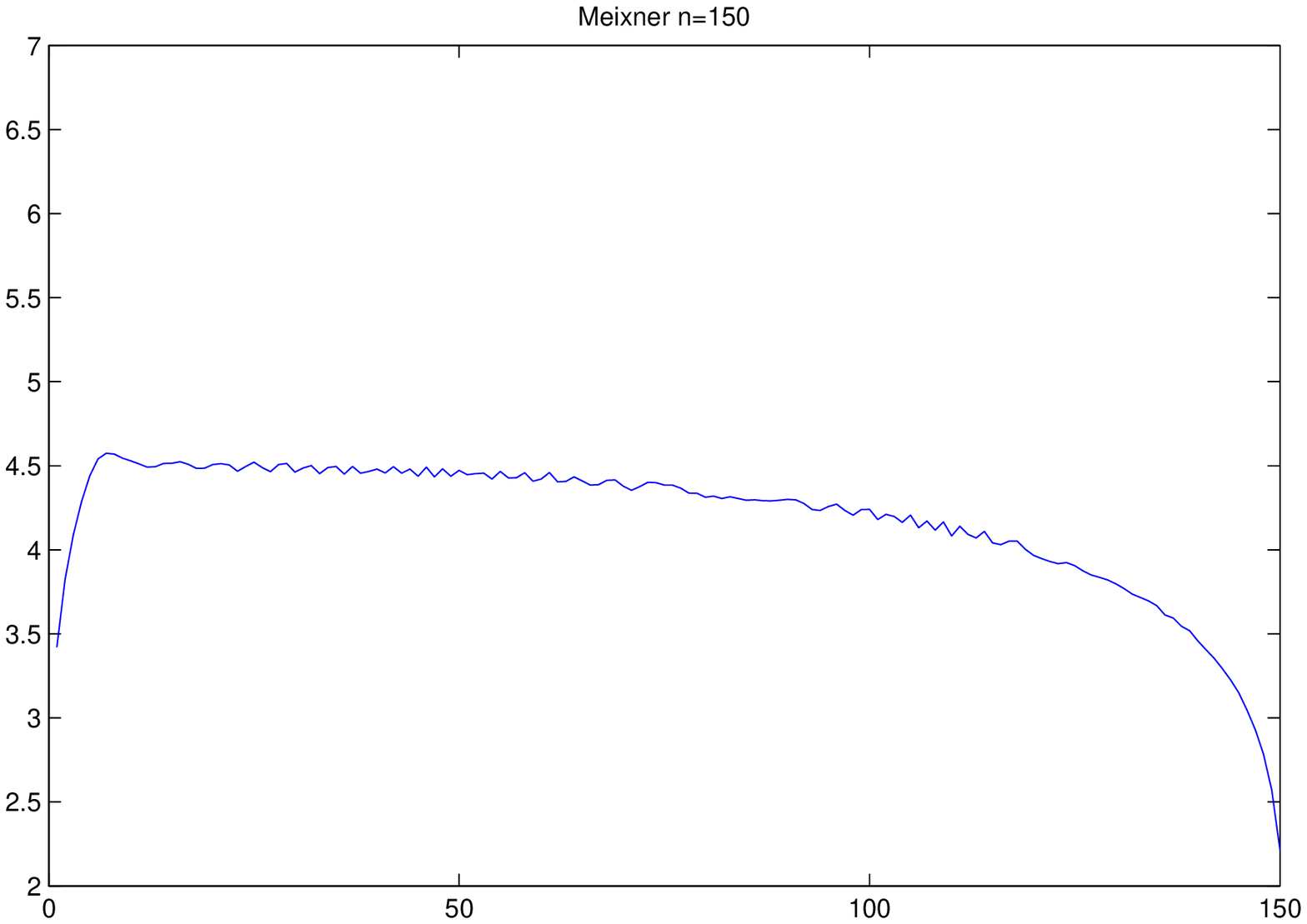}} &
\hspace{-1.2cm}\mbox{\includegraphics[scale=0.4]{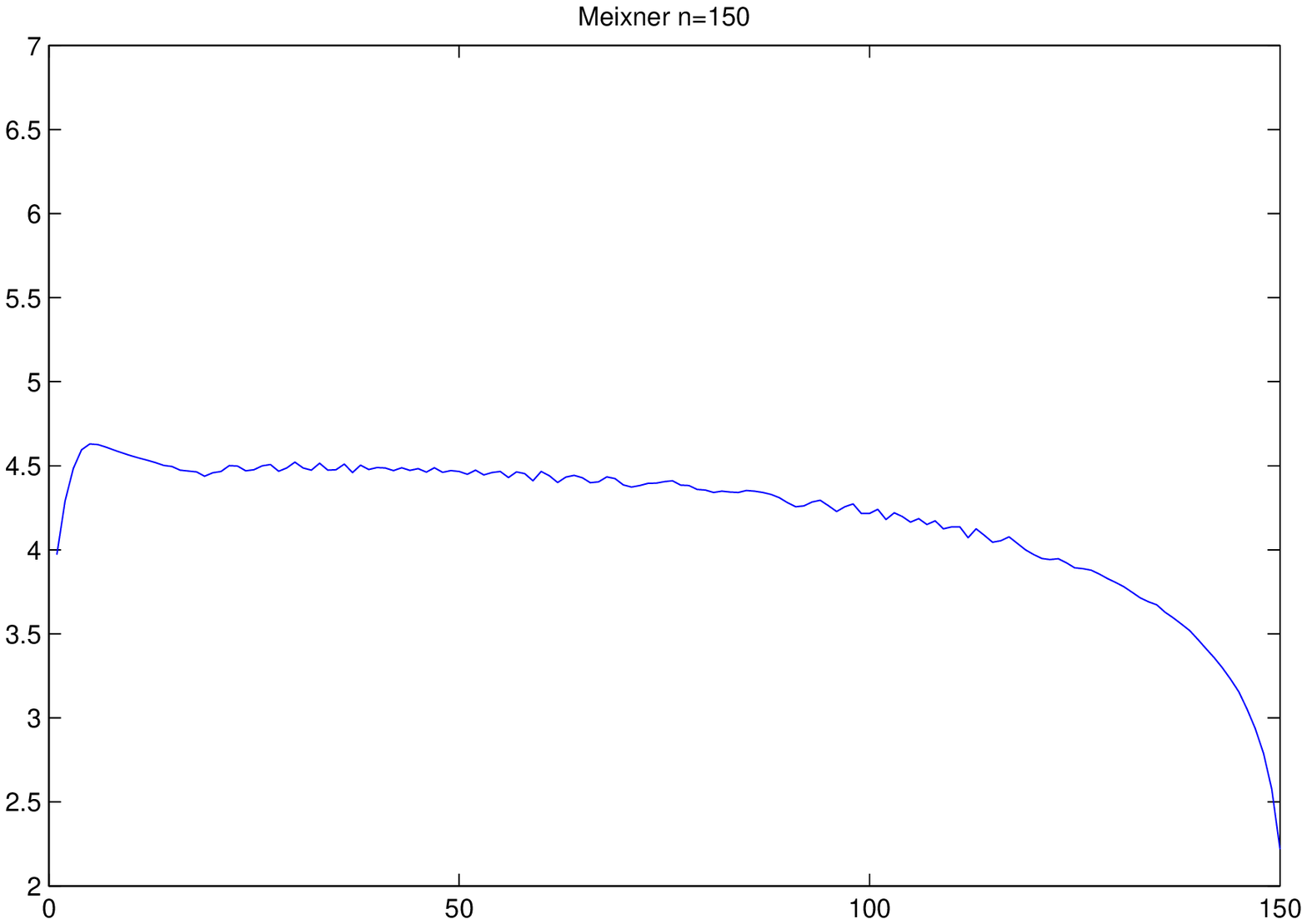}}
\end{tabular}
 \caption{Entropy $\S_{150,j}$  for Meixner polynomials $M_n^{(\beta ,c)} $ with $\beta  =3.4$, $c =0.2$  (left top), $\beta  =8.9$, $c =0.2$ (right top),
 $\beta  =3.4$, $c =0.8$ (left bottom), $\beta  =8.9$, $c =0.8$ (right bottom).
}\label{fig:Meixner}
\end{figure}
A qualitatively different behavior is observed for the Meixner polynomials $M_n^{(\beta ,c)}$, given by the recurrence relation \eqref{recurrence} with
$$
b_i = \frac{\left(i c ( i+\beta -1  )
\right)^{1/2}}{{ 1-c  }}\,, \quad a_i = \frac{(i-1)(1+c)+c\beta }{1-c}\,.
$$
Recall that they are orthogonal with respect to the discrete measure (see e.g.\ \cite[Chapter 6]{Ismail05})
$$
\mu = (1-c)^\beta \sum _{k=0}^\infty  \frac{(\beta )_k c^k}{ k!}   \,\delta _k\,.
$$
From Figure \ref{fig:Meixner} we observe that the value of the parameter $c $ has greater impact on the behavior of the entropy $\S_{n,j}$ in comparison with the parameter $\beta  $.

Finally, the evidence provided by all numerical experiments is sufficiently strong to conjecture that, after an appropriate rescaling and normalization, entropies $\S_{n,j}$ have a ``semiclassical'' limit as $n\to \infty$. The analysis of this asymptotic behavior is matter of a further research.

\section*{Acknowledgements}

AIA was partially supported by Programm N1 of DMSRAS and grant RFBR 08-01-00179 of Russian Federation.
JSD and RY were partially supported by Ministerio de Educaci\'{o}n y Ciencia, grant FIS 2005-00973, and by Junta de Andaluc\'{\i}a, grant FQM-207.
AMF acknowledges support from Ministerio de Educaci\'{o}n y Ciencia under grant
MTM2005-09648-C02-01, and from Junta de Andaluc\'{\i}a, grant
FQM-229. Additionally, JSD, AMF and RY were partially supported by the excellence grants FQM-481, and P06-FQM-01738 from Junta de Andaluc\'{\i}a.

AIA also wishes to acknowledge the hospitality of the Universities of Almer\'{\i}a and Granada, where this work was started.

%

\def\cprime{$'$}

\bigskip

\obeylines
\texttt{
A.I.\ Aptekarev (aptekaa@spp.keldysh.ru)
Keldysh Institute of Applied Mathematics,
Miusskaya Pl.\ 4, 125047 Moscow, RUSSIA
\medskip
J.S.\ Dehesa (dehesa@ugr.es)
Institute Carlos I of Theoretical and Computational Physics,
Granada University,
Campus de Fuentenueva, 18071 Granada, SPAIN
\medskip
A. Mart\'{\i}nez-Finkelshtein (andrei@ual.es)
Department of Statistics and Applied Mathematics, University of Almer\'{\i}a,
04120 Almer\'{\i}a, SPAIN, and
Institute Carlos I of Theoretical and Computational Physics,
Granada University, SPAIN
\medskip
R.\ Y\'{a}\~{n}ez (ryanez@ugr.es)
Institute Carlos I of Theoretical and Computational Physics,
Granada University,
Campus de Fuentenueva, 18071 Granada, SPAIN
}

\end{document}